\documentclass{article}
\usepackage{amsmath, amsthm, amssymb, amsfonts, textcomp, epsfig}
\newcommand{\mb}{\mathbb}
\newcommand{\zsp}{\mathbb{Z}/p\,\mb{Z}}
\newcommand{\zspx}{\mathbb{Z}/p\,\mb{Z}[X]}
\newcommand{\tk}{t_k^{(0)}}
\newcommand{\tku}{t_k^{(1)}}
\newcommand{\qq}{\qquad\text{mod}\,p^2}

\newcommand{\dk}{\delta_0(k)}
\newcommand{\di}{\delta_0(i)}
\newcommand{\djz}{\delta_0(j)}
\newcommand{\dku}{\delta_1(k)}
\newcommand{\diu}{\delta_1(i)}

\newcommand{\s}{\sum_{k=1}^{p-1}}
\newcommand{\si}{\sum_{i=1}^{p-1}}
\newcommand{\mdt}{\;\;\text{mod}\,p^3}
\newcommand{\sti}{\left[\begin{array}{l}\;\;\;\,p\\2n+1\end{array}\right]}
\newcommand{\stp}{\left[\begin{array}{l}\;\;\;\,p\\2n+2\end{array}\right]}

\newcommand{\bn}{\binom}
\author{Claire Levaillant}
\title{Wilson's theorem modulo $p^2$ derived from Faulhaber polynomials}
\begin{document}
\maketitle

\noindent Abstract. First, we present a new proof of Glaisher's formula dating from $1900$ and concerning Wilson's theorem modulo $p^2$. Our proof uses $p$-adic numbers and Faulhaber's formula for the sums of powers ($17$th century), as well as more recent results on Faulhaber's coefficients obtained by Gessel and Viennot. Second, by using our method, we find a simpler proof than Sun's proof regarding a formula for $(p-1)!$ modulo $p^3$, and one that can be generalized to higher powers of $p$. Third, we can derive from our method a way to compute the Stirling numbers $\left[\begin{array}{l}p\\s\end{array}\right]$ modulo $p^3$, thus improving Glaisher and Sun's own results from $120$ years ago and $20$ years ago respectively. Last, our method allows to find new congruences on convolutions of divided Bernoulli numbers and convolutions of divided Bernoulli numbers with Bernoulli numbers.

\section{Introduction}
The paper arose from an analogy between two polynomials, namely
\begin{center}
$f(X)=X^{p-1}-1\in\mathbb{Z}/p\,\mathbb{Z}[X]$ and $g(X)=X^{p-1}+(p-1)!\in\mb{Z}_p[X],$
\end{center}
where $p$ is a prime number, $\mb{Z}/p\,\mb{Z}$ denotes the field with $p$ elements and $\mb{Z}_p$ denotes the ring of $p$-adic integers, see \cite{GO}.\\\\
The analogy is concerned with the way each polynomial factors and with the congruences properties that can be derived in each case from the relations between the coefficients and the roots. \\\\Below, we describe the analogy in details.
\begin{itemize}
\item First we look at the factorization of $f(X)$ and the properties modulo $p$ that can be derived from it.
By Fermat's little theorem we have,
$$k^{p-1}=1\;\;\text{in}\;\;\mb{Z}/p\,\mb{Z}$$
for all integers $k$ with $1\leq k\leq p-1$. It provides the $p-1$ roots of $f$, and so $f$ factors in $\zspx$ as
$$f(X)=(X-1)(X-2)\dots\,(X-(p-1))$$
Let us introduce some notations concerning the divided factorial which shall be useful throughout the paper.
\newtheorem{Notation}{Notation}
\begin{Notation}
Let $r$ be an integer with $1\leq r\leq\,p-2$ and let $i_1,\,i_2,\dots,\,i_r$ be $r$ integers with $1\leq i_k\leq\,p-1$ for all $k$ with $1\leq k\leq r$. \\\\We define the divided factorial $(p-1)!^{i_1,\dots,\,i_r}$ as
$$(p-1)!^{i_1,\dots,\,i_r}:=\frac{(p-1)!}{\prod_{k=1}^r i_k}$$
\end{Notation}
From looking at the constant coefficient of $f$ in both factored and expanded forms, we retrieve Wilson's theorem which reads
$$(p-1)!=-1\;mod\;p$$
This constitutes one of many proofs for Wilson's theorem.\\
By looking at the other coefficients in both factored and expanded forms, we derive more relations, namely,
\begin{eqnarray*}
\sum_{k=1}^{p-1}(p-1)!^k&=&0\qquad mod\,p\qquad\text{(Coefficient in $X$)}\\
\sum_{1\leq i<j\leq p-1}(p-1)!^{i,j}&=&0\qquad mod\,p\qquad\text{(Coefficient in $X^2$)}\\
&\vdots&\qquad\qquad\qquad\qquad\qquad\vdots\\
\sum_{1\leq i_1<i_2<\dots<i_r\leq p-1}(p-1)!^{i_1,\dots,\,i_r}&=&0\qquad mod\,p\qquad\text{(Coefficient in $X^r$)}\\
&\vdots&\qquad\qquad\qquad\qquad\qquad\vdots\\
\sum_{1\leq i_1<i_2<\dots<i_{p-2}\leq p-1}(p-1)!^{i_1,\dots,\,i_{p-2}}&=&0\qquad mod\,p\qquad\text{(Coefficient in $X^{p-2}$)}
\end{eqnarray*}

The numbers to the left hand side are also the respective unsigned coefficients of $x^2,x^3,\,\dots,\,x^{r+1},\,\dots,\,x^{p-1}$ in the falling factorial
$$x(x-1)(x-2)\dots(x-(p-1)),$$
that is the unsigned Stirling numbers of the first kind.
$$\left[\begin{array}{l}p\\2\end{array}\right],\,\left[\begin{array}{l}p\\3\end{array}\right],\dots,\left[\begin{array}{l}\;\;\;p\\r+1\end{array}\right],\dots,\left[\begin{array}{l}\;\;\;p\\p-1\end{array}\right]$$
We summarize the preceding equalities by
$$\forall 2\leq k\leq p-1,\;p\;|\left[\begin{array}{l}p\\k\end{array}\right]$$
These numbers were introduced by James Stirling in the $18$th century. The bracket notation gets promoted by Donald Knuth by analogy with binomial coefficients. Namely, we have the well-known divisibility relations
$$\forall 1\leq k\leq p-1,\;p\;|\left(\begin{array}{l}p\\k\end{array}\right),$$
a straightforward consequence of the definition of the binomial coefficients and of the Gauss lemma. \\
Note that another equivalent definition for the unsigned Stirling numbers of the first kind is the following:
$\left[\begin{array}{l}n\\s\end{array}\right]$ counts the number of permutations of $n$ elements that decompose into a product of $s$ disjoint cycles, see \cite{LEV} for a quantum combinatorial proof. \end{itemize}
While there exist lots of proofs for Wilson's theorem, little is advertised about Wilson's theorem modulo $p^2$ or higher powers. It is even said at the end of \cite{GR}: "There is no explicit formula for $(p-1)!\,\text{mod}\,p^2$, ie the integer $n_1$ in $(p-1)!\equiv\,-1+n_1p\;\text{mod}\,p^2$ is no evident function of $p$". \\In \cite{HW}, Wilson's theorem modulo $p^2$ is only stated as
$$(p-1)!\equiv (-1)^{\frac{p-1}{2}}\,2^{2p-2}\,\bigg(\frac{p-1}{2}!\bigg)^2\qquad(\text{mod}\,p^2)$$ However, British mathematician J.W.L. Glaisher, son of meteorologist and pioneer both of weather forecasting and of photography James Glaisher who in particular held the world record for the highest altitude reached in a balloon, found a formula for $(p-1)!\,\text{mod}\,p^2$ in $1900$ \cite{GL} that uses Bernoulli numbers. His method is elegant, smart and straightforward, but does not generalize to higher powers of $p$.\\
A hundred years later, in $2000$, Z-H Sun provides a different perspective, see \cite{SU2}. His method even allows him to compute $(p-1)!\,\text{mod}\,p^3$, thus generalizing Glaisher's result whose work had only led to a formula for $(p-1)!\,\text{mod}\,p^2$. Sun's result modulo $p^3$ is expressed in terms of divided Bernoulli numbers. Sun's machinery is heavy but really beautiful and well exposed. First, he determines the generalized harmonic numbers $H_{p-1,k}$ modulo $p^3$. For that, he uses a battery of tools, namely: Euler's theorem ; Bernoulli's formula for the sums of powers modulo $p^3$ together with von Staudt--Clausen theorem \cite{VS}\cite{CL} that determines the fractional part of Bernoulli numbers ; K\"ummer's congruences from $1851$ involving Bernoulli numbers \cite{KU} and their generalizations modulo $p^2$ by Sun himself \cite{SU2} in the case when $k\leq p-3$, and an unpublished result of \cite{SU3} giving formulas for $p\,B_{k(p-1)}$ modulo $p^2$ and $p^3$ respectively for dealing with $k=p-2$ and $k=p-1$ respectively. Second, he uses Newton's formulas together with Bernoulli's formula modulo $p^2$ in order to derive the Stirling numbers modulo $p^2$. From there he obtains in particular a congruence for $(p-1)!\,\text{mod}\,p^3$ that was first proven by Carlitz \cite{CAR}. Likewise, by using Newton's formulas with generalized harmonic numbers and what Sun denotes as the conjugates of the Stirling numbers, he finds a formula for the conjugate Stirling numbers modulo $p^2$ followed by a "conjugate Carlitz congruence". By combining the various identities, he then derives a pioneering formula for $(p-1)!\,\text{mod}\,p^3$ that is expressed only in terms of Bernoulli numbers.\\
Yet another twenty years later, we offer a third perspective which does not make any use of the Newton formulas. The polynomial $g$ with $p$-adic integer coefficients introduced earlier happens to be an even more powerful tool than the polynomial $f$ for better prime precision. It highlights once again the beauty and efficiency of $p$-adic numbers when investigating local properties at powers of primes. \\

Explicitly:
\begin{itemize}
\item We now investigate the factorization of $g(X)$. \\Let $k$ be an integer with $1\leq k\leq\,p-1$.\\
First and foremost, we have
\begin{equation}g(k)=k^{p-1}+(p-1)!\in\,p\,\mb{Z}_p\end{equation}
by definition of $g$ and by using Wilson's theorem modulo $p$.\\
Moreover, we have
\begin{equation}
g^{'}(k)=(p-1)\,k^{p-2}\not\in\,p\,\mb{Z}_p
\end{equation}
as $k^{p-2}=k^{-1}\neq\,0$ in $\zsp$.\\
By Hensel's lemma \cite{GO}, each nonzero $k$ of $\zsp$ lifts to a unique root $x_k$ of $g$ such that
$$k-x_k\in\,p\,\mb{Z}_p$$
Write $$x_k=k+p\,t_k\qquad\text{with $t_k\in\mb{Z}_p$}$$
Then $g$ factors in $\mb{Z}_p[X]$ as follows.
$$g(X)=(X-1-p\,t_1)\dots\,(X-(p-1)-p\,t_{p-1})$$
We will see that by looking at the constant coefficient of $g$ modulo $p^2\,\mb{Z}_p$, we obtain the next coefficient in the $p$-adic expansion of $(p-1)!$, namely
\newtheorem{Theorem}{Theorem}
\begin{Theorem}
$$(p-1)!=-1+p\;\sum_{k=1}^{p-1}\delta_0(k)\qquad\text{mod}\;p^2,$$ with $\delta_0(k)$ defined for each $1\leq k\leq p-1$ by
$$k^{p-1}=1+p\;\delta_0(k)\qquad\text{mod}\,p^2$$
\end{Theorem}
Working out the sum of Theorem $1$ modulo $p^2$ by using the Faulhaber polynomials for the sums of powers and by using Faulhaber's neat statement on the relationship between the two trailing coefficients of the polynomial for odd powers (work done in the $17$th century, \cite{JF}), we obtain Wilson's theorem one step further, namely modulo $p^2$.
\begin{Theorem}
Let $p$ be an odd prime number. Set $p=2l+1$ and $a=\frac{p(p-1)}{2}$. \\
Let $c_1(l)$ be the trailing coefficient in the Faulhaber polynomial
\begin{equation}\sum_{k=1}^{p-1}\,k^p\,=c_l(l)\,a^{l+1}+\dots\,+c_2(l)\,a^3+c_1(l)\,a^2\end{equation}
We have,
$$\begin{array}{l}(i)\;\;\forall\,i\in\,\lbrace\,1,\dots,\,l\rbrace,\;c_i(l)\in\mb{Z}_p\\\\
(ii)\;\,(p-1)!=\frac{1}{2}\,c_1(l)-p\;\;\;\text{mod}\,p^2\end{array}$$
\end{Theorem}
\newtheorem{Corollary}{Corollary}
\begin{Corollary} (Glaisher $1900$ \cite{GL})
$$(p-1)!=p\,B_{p-1}-p\qquad\text{mod}\,p^2$$
\end{Corollary}
Corollary $1$ can be derived directly from Theorem $1$ by using the Bernoulli formula for the sum of powers smartly modulo $p^2$ (like Sun did in \cite{SU2}, see the proof of his congruence $(5.1)$) together with von Staudt--Clausen theorem. \\For future reference, Bernoulli's formula for the sums of powers reads:
$$\sum_{k=1}^pk^m=\frac{1}{m+1}\sum_{j=0}^m\binom{m+1}{j}B_j\,p^{m+1-j},$$
where we wrote the formula using Bernoulli numbers of the second type, that is $B_1=\frac{1}{2}$ instead of $B_1=-\frac{1}{2}$.\\
As for von Staudt--Clausen's theorem it gets described just a little further below.

Corollary $1$ can also be derived from Theorem $2$. Indeed, the Faulhaber coefficients of $(3)$ are related to the coefficients $A_k^{(l+1)}$ of \cite{KN}, which can themselves be computed in terms of Bernoulli numbers by the formula first proven by Gessel and Viennot (see \cite{GV} and \cite{KN}). Edwards noticed for fixed $m$ that the coefficients $A_k^{(m)}$ can be defined recursively, thus are obtainable by inverting a lower triangular matrix, see \cite{Ed}. From there, Gessel and Viennot found a neat combinatorial interpretation for these coefficients. From their magnificent result, it can then be conveniently derived that some of the rational Faulhaber coefficients are $p$-adic integers (so are for instance the Faulhaber coefficients appearing in our Theorem $2$), a useful data to have when working modulo powers of $p$.\\
We will see and point out from now on that by using a result of Derby \cite{NI}, Glaisher's formula provides an efficient way of computing $(p-1)!$ modulo $p^2$.\\
Finally, we note that Glaisher's formula is consistent with the fact that \begin{equation}p\,B_{p-1}=-1\;\text{mod}\,p\,\mb{Z}_p\end{equation} The latter fact is indeed a consequence of the previously mentioned von Staudt--Clausen theorem. This theorem dating from $1840$ and independently proven by von Staudt and by Clausen asserts that the Bernoulli numbers $B_{2k}$ sum to zero when added all the fractions $\frac{1}{q}$ with $q$ prime such that $q-1$ divides $2k$ (recall that $B_{2k+1}=0$ when $k\geq 1$). Von Staudt, German mathematician, was a student of Gauss at G\"ottingen. Danish mathematician Clausen, while illiterate at age $12$ and astronomer by background, was so talented that Gauss once said when Clausen was in life trouble:
\begin{center}\textit{"Es w\"are doch sehr zu beklagen, wenn sein wirklich ausgezeichnetes Talent f\"ur abstracte Mathematik in der Verk\"ummerung so ganz zu Grunde ginge"}\end{center}
Amongst his many achievements, Clausen calculated in $1847$ the first $247$ decimals of $\pi$ and in $1854$ he decomposed the Fermat number $F_6$ into prime numbers, thus uncovering a prime number that was the biggest of all primes formerly known, namely $67280421310721$. \\
In mathematical words, von Staudt--Clausen's theorem reads
$$B_{p-1}+\sum_{\begin{array}{l}\;\;q\in\mathcal{P}\\q-1|p-1\end{array}}\frac{1}{q}=0\;\;\text{mod}\,1$$
Then, (multiplying by $p$),
$$p\,B_{p-1}+1+p\sum_{\begin{array}{l}\;\;q\in\mathcal{P}\\\;\;q\neq\,p\\q-1|p-1\end{array}}\frac{1}{q}=0\;\;\text{mod}\,p\,\mb{Z}_p$$
Leading to the result stated in $(4)$.
Note that von Staudt--Clausen theorem specifically applied to the Bernoulli number $B_{p-1}$ implies that the denominator of $B_{p-1}$ consists of $p$ times other primes $q$ such that $q-1|p-1$.
The Agoh-Giuga conjecture, named after independent Japanese mathematician Takashi Agoh and Italian mathematician Giuseppe Giuga asserts that $(4)$ holds only for prime numbers.
\newtheorem{Conjecture}{Conjecture}
\begin{Conjecture} (Agoh-Giuga \cite{AG}\cite{GI})\\
$n$ is prime if and only if $n\,B_{n-1}=-1\;\text{mod}\,n$
\end{Conjecture}
If Conjecture $1$ holds, we must show that $n\,B_{n-1}\neq -1\;\text{mod}\,n$ for $n$ composite in order to prove it. Given $n$ composite, if $n\,B_{n-1}=0\;\text{mod}\,n$, then we are done. Otherwise, from von Staudt-Clausen's theorem, there must exist a prime divisor $p$ of $n$ such that $p-1|n-1$. We say that $n$ is "Carmichael at $p$". A Carmichael number is a composite square-free integer that is Carmichael at all its prime divisors. Equivalently,
a Carmichael number is a composite number which is a Fermat pseudo-prime to every base \cite{CARM}, that is
$$\forall a\in\mb{Z},\;a^n=a\;\text{mod}\,n$$
In other words, the Carmichael numbers are the "Fermat liars". In particular, a Carmichael number must be odd. The equivalence between both definitions was shown in $1899$ by Korselt \cite{KO}.
Giuga has linked the property of Conjecture $1$ to Carmichael numbers and Giuga numbers.
\newtheorem{Definition}{Definition}
\begin{Definition}
$n$ is a Giuga number if $$\sum_{\begin{array}{l}p\in\mathcal{P}\\p|n\end{array}}\frac{1}{p}-\prod_{\begin{array}{l}p\in\mathcal{P}\\p|n\end{array}}\frac{1}{p}\in\mb{N}$$
Equivalently, for each prime divisor $p$ of $n$, we have $p\,\big|\frac{n}{p}-1$.
\end{Definition}
Giuga has proven the following result.
\newtheorem*{Thm}{Theorem}
\begin{Thm} (Giuga 1950 \cite{GI})\\
A composite integer $n$ satisfies $nB_{n-1}=-1\;\text{mod}\;n$ if and only if \\$n$ is both a Giuga number and a Carmichael number.
\end{Thm}
\noindent In light of Giuga's theorem, the Agoh-Giuga conjecture asserts the non-existence of composite numbers that are both Giuga numbers and\\ Carmichael numbers. We will see in Section $3$ of the paper that an odd composite Giuga number such that for any prime divisor $p$ of $n=p\,m$, $p$ does not divide $B_{m-1}$ would provide a counter-example to the Agoh-Giuga conjecture. Indeed, we have the following fact.
\newtheorem{Fact}{Fact}
\begin{Fact}
Let $n$ be an odd composite Giuga number. \\Let $p\in\mathcal{P}$ with $p|n$. Write $n=pm$, some integer $m$.\\
Then, $$p-1|m-1\;\text{or}\;p\big|B_{m-1}$$
\end{Fact}
In Section $3$, we prove Fact $1$ without reference to Adams'theorem \cite{AD}\cite{JO}.\\

\indent Going now back to the polynomial $g$ and looking at the respective coefficients in $X,\,X^2,\dots,\,X^{p-2}$ modulo $p^2$ in both factored and expanded forms of $g$ (just like we did with polynomial $f$ modulo $p$ only), we are able to express the Stirling numbers one prime power further, in terms of a difference between a sum of powers and a generalized harmonic number, like gathered in the following theorem.
\begin{Theorem} Let $S_{p-1,p-k}$ denote the sum of the $(p-k)$-th powers of the first $(p-1)$ integers, that is
$$S_{p-1,p-k}:=\sum_{l=1}^{p-1}l^{\,p-k}$$ Let $H_{p-1,k-1}$ denote the generalized harmonic number, using a standard notation. \\We have,
$$\forall\,2\leq k\leq p-1,\,\left[\begin{array}{l}p\\k\end{array}\right]= S_{p-1,p-k}-\,H_{p-1,k-1}\;\;\text{mod}\,p^2$$
\end{Theorem}
The proof of Theorem $3$ uses a result by Bayat which is a generalization of Wolstenholme's theorem. Wolstenholme's theorem asserts that
$$H_{p-1,1}=0\;\text{mod}\,p^2\;\; \text{and}\;\;H_{p-1,2}=0\;\text{mod}\,p$$
Back in $1862$, Wolstenholme proved the result \cite{WO} (which is since then referred to as Wolstenholme's theorem) as an intermediate result for showing that
$$\binom{2p-1}{p-1}\equiv 1\;\text{mod}\,p^3$$
The latter result is also commonly referred to as Wolstenholme's theorem. \\
A $150$ years later, generalizing Wolstenholme's theorem was still of interest to mathematicians. In $1997$, Bayat generalized the part of Wolstenholme's theorem related to the generalized harmonic numbers. He showed the following fact.
\newtheorem{Result}{Result}
\begin{Result} (Bayat, 1997 \cite{BA})
Let $m$ be a positive integer and let $p$ be a prime number with $p\geq m+3$. Then,
$$\sum_{k=1}^{p-1}\frac{1}{k^m}\;\equiv\;\begin{cases}0\;\;\text{mod}\;p\;\;\;\text{if $m$ is even}\\
0\;\;\text{mod}\;p^2\;\;\;\text{if $m$ is odd}
\end{cases}$$
\end{Result}
Note, and this becomes important in our discussion later on, in the proof of Theorem $3$, we only need Bayat's result modulo $p$. \\
From Theorem $3$ and using this time a more recent generalization of Wolstenholme's theorem by Sun, we then derive divisibility properties or congruences modulo $p^2$ for the\\$\begin{array}{l}\end{array}$\\unsigned Stirling numbers of the first kind $\left[\begin{array}{l}p\\2\end{array}\right],\left[\begin{array}{l}p\\3\end{array}\right],\dots,\,\left[\begin{array}{l}\;p\\p-1\end{array}\right]$. \\

The results get summarized in the corollary below.
\begin{Corollary}
Let $n$ be an integer and $p$ be a prime. We have
$$\begin{array}{l}
\left[\begin{array}{l}\;\;\;p\\2n+2\end{array}\right]=\begin{cases}
0\;\;\,\;\,\,\;\text{mod}\;p^2 & \text{if}\;\;2n+1\leq p-4\\
-\frac{p}{2}\;\;\text{mod}\;p^2& \text{if}\;\;2n+1=p-2\\
\end{cases}\\\\
\left[\begin{array}{l}\;\;\;p\\2n+1\end{array}\right]=\begin{cases}
\frac{p}{2n+1}\;B_{p-2n-1}\;\text{mod}\;p^2&\text{if}\;\;p>2n+1\;\,\text{and}\;\,n\geq 1\\
1\qquad\qquad\qquad\!\text{mod}\;p^2&\text{if}\;\;p=2n+1\\
p\,B_{p-1}-p\qquad\!\text{mod}\;p^2&\text{if}\;\;n=0
\end{cases}
\end{array}$$
\end{Corollary}
The result of Corollary $2$ also appears in the original and subtle work of Glaisher from $1900$. It is stated as follows using Glaisher's notations.
\begin{Result} (Glaisher, 1900 \cite{GL})\\
Denote by $A_r$ the sum of products of r integers among the first $(p-1)$ integers. We have,
\begin{eqnarray}
\frac{A_1}{p}&=&-\frac{1}{2}\;\;\qquad\text{mod}\,p\\
\frac{A_{2k+1}}{p}&=&\;\;\;0\;\;\qquad\text{mod}\,p\\
\frac{A_{2k}}{p}&=&-\frac{V_{2k}}{2k}\;\;\;\;\text{mod}\,p
\end{eqnarray}
\end{Result}
The link between Glaisher's notations and our notations is
$$\left[\begin{array}{l}p\\k\end{array}\right]=A_{p-k}$$
As for the Bernoulli numbers, his V-notation corresponds to our B-notation.


Contrary to what happens with our own method for finding Wilson's theorem modulo $p^2$, the formulas providing the Stirling numbers modulo $p^2$ must be used within Glaisher's method in order to derive Wilson's theorem modulo $p^2$. Thus, our method for proving Wilson's theorem modulo $p^2$ appears somewhat more straightforward than Glaisher's way since our way uses only p-adic numbers (cf Theorem $1$) and the Bernoulli's formula modulo $p^2$ for the sums of powers.

Finally, we note that the result of Corollary $2$ also appears in Sun's work, based on Newton's formulas and based also on divisibility properties at prime $p$ of the sums of powers. Still using Glaisher's notations and writing simply $S_k$ for $S_{p-1,k}$ (thus lightening our own notation from Theorem $3$), Newton's formulas relate $A_k$ and $S_k$ for each $k=1,\dots,p-1$ by (cf \cite{JAC})
\begin{equation}A_k=\frac{(-1)^{k-1}}{k}\bigg(S_k+\sum_{r=1}^{k-1}(-1)^r\,A_r\,S_{k-r}\bigg)\end{equation}
As we know, $p$ divides $A_r$ for all $1\leq r\leq\,k-1$ when $k\leq p-1$. Moreover, it is clear from Faulhaber's formula that $p$ (and even $p^2$) divides $S_3,\dots,\,S_{p-2}$ for the odd indices once we know that the Faulhaber coefficients that are involved are $p$-adic integers. And of course $p$ divides $S_1$ by a direct calculation. Dealing next with the even indices,
from Faulhaber's formula applied with even powers, namely
\begin{equation}
\sum_{k=1}^{p-1}k^{2l}=\frac{p-1+\frac{1}{2}}{2l+1}\;\,\big(2c_1(l)a+3c_2(l)a^2+\dots+(l+1)c_l(l)a^l\big)
\end{equation}
we can show by using the work of Gessel and Viennot previously evoked that the Faulhaber coefficients that are present in Formula $(9)$ are all $p$-adic integers. Thus, for each integer $l$ with $1\leq l\leq \frac{p-3}{2}$, we have
\begin{equation}
S_{2l}=\frac{1}{2}\,\frac{p\,c_1(l)}{2l+1}\;\text{mod}\,p^2
\end{equation}
In particular, it follows that $p$ divides $S_{2l}$ for such $l$'s.
Gathering these divisibility facts, we may now assert that for each $k=1,\dots,\,p-1$, the sum in $(8)$ is congruent to $0$ modulo $p^2$. And so Sun obtains
\begin{equation}
A_k=\frac{(-1)^{k-1}}{k}\;S_k\;\;\text{mod}\,p^2\;\;\;\forall\,k=1,\dots,\,p-1
\end{equation}
Combining our Theorem $3$ and Sun's result in $(11)$, we can now relate the generalized harmonic numbers to the sums of powers modulo $p^2$ in the following way.
\begin{Corollary}
\begin{equation}
H_{p-1,p-k-1}=\bigg(1+\frac{(-1)^{k}}{k}\bigg)\,S_k\;\;\;\text{mod}\;\,p^2\;\;\;\forall\,k=1,\dots,p-2
\end{equation}
\end{Corollary}
Setting $k=1$, this shows in particular that \begin{equation}H_{p-1,p-2}=0\;\;\text{mod}\,p^2\end{equation}
More generally, since we have seen that $p^2$ divides $S_3,\dots,\,S_{p-2}$, we also have
\begin{Corollary}
\begin{equation}
H_{p-1,k}=0\;\text{mod}\,p^2\qquad\forall k=1,3,\dots,p-4,p-2\;\text{with}\;k\;\text{odd}
\end{equation}
\end{Corollary}
The latter result to the exception of $(13)$ can be found in \cite{HW} on page $103$ with a proof unrelated to ours.
The case $k=1$ is usually referred to as Wolstenholme's theorem.

Here, we have thus recovered Bayat's result for the modulus $p^2$ case if we gather $(10)$, $(12)$ and $(14)$. Note with $(13)$ that our result is an improvement to Bayat's result since Bayat does not deal with the generalized harmonic number $H_{p-1,p-2}$. More recently in $2011$, Romeo Mestrovic proved in \cite{RO} that
\begin{equation}\binom{2p-1}{p-1}\equiv 1-2p\,H_{p-1,1}+4p^2\;A_2^{*}\;\;(mod\;p^7),\end{equation}
where $A_2^{*}$ denotes the conjugate of $A_2$ using Sun's notations. See also \cite{RO2} for more complements on the topic.
James P. Jones has conjectured that the converse of Wolstenholme's theorem holds and there is past and ongoing research on the matter, see for instance in chronological order \cite{MC}, \cite{HT}, \cite{HU}.

Closing this digression and going back to congruence $(10)$, it rewrites in terms of Bernoulli numbers as
$$S_{2l}=p\,B_{2l}\;\text{mod}\,p^2$$ 
see forthcoming $\S\,2$. Using this rewriting and the fact that $B_j=0$ for odd $j>1$, we derive further from $(12)$ the general formula
\begin{Corollary}
\begin{equation}
H_{p-1,p-k-1}=\bigg(1+\frac{(-1)^{k}}{k}\bigg)\;p\,B_k\;\;\;\text{mod}\;\,p^2\;\;\;\forall\,k=1,\dots,p-2
\end{equation}
\end{Corollary}
Note the formula still holds for $k=1$ though $B_1\neq 0$.\\

Recall that when $k\leq p-2$, the denominator of $B_k$ is not divisible by $p$ since $p-1$ does not divide $k$.
Then, modulo $p^2$, the right hand side of Equality $(16)$ is precisely
$$\frac{p-k-1}{p-k}\,p\,B_k,$$
which is nothing else than the formula provided in Corollary $5.1$ of Sun in \cite{SU2}.
This is a result in Sun, which though generalizable to $p^3$, had required a fair amount of work. However, we do not reveal any congruence modulo $p^2$ for $H_{p-1,p-1}$ by our method, whereas Sun does deal with that case successfully. And like three mathematicians working in mirrors, Glaisher himself has worked out the generalized harmonic numbers to the modulus $p^2$ and by times $p^3$. He noticed from
\begin{equation}
(X-1)(X-2)\dots(X-(p-1))=X^{p-1}-A_1\,X^{p-2}+\dots\,-A_{p-2}X+A_{p-1}
\end{equation}
that $$1,\frac{1}{2},\frac{1}{3},\dots,\frac{1}{p-1}$$ are the roots of the polynomial
$$A_{p-1}X^{p-1}-A_{p-2}X^{p-2}+\dots+A_2X^2-A_1X+1$$
From there, by using the Newton formulas and some divisibility properties of the Stirling numbers, namely $$\left\lbrace\begin{array}{l}p|\left[\begin{array}{l}p\\k\end{array}\right]\;\text{for every integer}\;k\;\text{with}\;2\leq k\leq p-1\\\\ \negthickspace p^2|\left[\begin{array}{l}p\\2l\end{array}\right]\;\text{for every integer}\;l\;\text{with}\;1\leq l\leq\frac{p-3}{2}\end{array}\right.,$$ he derives as part of his work, where his $H_k$ is our $H_{p-1,k}$,
\begin{Theorem} (Glaisher $1900$ \cite{GL})
$$\begin{array}{ccccc}
H_1&\equiv&-A_{p-2}&\text{mod}&p^3\\
H_2&\equiv&2A_{p-3}&\text{mod}&p^2\\
H_3&\equiv&-3A_{p-4}&\text{mod}&p^3\\
&\vdots &&&\\
&\vdots &&&\\
H_{p-3}&\equiv&(p-3)A_2&\text{mod}&p^2
\end{array}$$
with the modulus being alternately $p^3$ and $p^2$, and
$$\begin{array}{ccccc}
H_{p-2}&\equiv&-p-(J-\frac{3}{2})\;p^2&\text{mod}&p^3\\
H_{p-1}&\equiv&-1-(J-1)\;p&\text{mod}&p^2\end{array},$$
with the value of J given by
\begin{equation}J=-1+V_{p-1}+\frac{1}{p},\end{equation}
that is $J$ is the second coefficient in the p-adic expansion of $(p-1)!$. Namely, we have
\begin{equation}(p-1)!=-1+p\,J\;\;\text{mod}\,\;p^2\end{equation}
\end{Theorem}

Thus, our work for finding the Stirling numbers with (with respect to our notations) even indices modulo $p^2$ can be deduced from Theorem $3$ joint with Hardy and Wright's congruence for generalized harmonic numbers with odd indices modulo $p^2$ (or Bayat's result modulo $p^2$). Once this is achieved, finding the Stirling numbers with (with respect to our notations) odd indices modulo $p^2$ follows jointly from Theorem $3$ together with Glaisher's formulas for generalized harmonic numbers modulo $p^2$, without reference to Sun's later work. \\Glaisher actually derives explicit formulas for the generalized harmonic numbers $H_k$ modulo $p^3$ for odd indices $k$ because as part of his unlimited and groundbreaking creativity, not only had he calculated the Stirling numbers modulo $p^2$ in terms of Bernoulli numbers, but he had also found a way to lift this calculation to the next modulus $p^3$ as far as the odd indices are concerned.
He had shown in his earlier paper \cite{GL1} the following fact.
\begin{Result} (Glaisher, 1900 \cite{GL1}) Let $r$ be an odd integer with $3\leq r\leq p-2$. We have,
\begin{equation}A_r\;\equiv\;\frac{p\,(p-r)}{2}\;A_{r-1}\;\;\text{mod}\,p^3\end{equation}
\end{Result}
When $r$ is odd, $r-1$ is even, whence by $(7)$,
\begin{equation}A_{r-1}=-p\;\frac{V_{r-1}}{r-1}\;\;\text{mod}\,p^2\end{equation}
Gathering $(20)$, $(21)$ and the fact that $A_{r-1}=0\;\text{mod}\,p$, it leads to
\begin{Result} (Glaisher, 1900 \cite{GL})
Let $r$ be an odd integer with $3\leq r\leq p-2$. We have,
\begin{equation}
A_r\;\equiv\;\frac{p^2\,r}{2(r-1)}\;V_{r-1}\;\;\text{mod}\,p^3
\end{equation}
\end{Result}
If $r$ is odd, so is $p-1-r$ and if $1\leq r\leq p-4$, then $3\leq p-1-r\leq p-2$, thus by $(22)$ and von Staudt-Clausen's theorem, we have
\begin{equation}
A_{p-1-r}\equiv\;\frac{p^2(r+1)}{2(r+2)}\;V_{p-r-2}\;\;\text{mod}\,p^3
\end{equation}
From Theorem $4$, we have $H_r=-r\,A_{p-r-1}\;\;\text{mod}\,p^3$ for each odd $r$ with $1\leq r\leq p-4$.
Hence,
$$H_r=-\frac{p^2\,r(r+1)}{2(r+2)}\;V_{p-r-2}\;\;\text{mod}\,p^3$$
We now state Glaisher's result in our own notations, emphasizing the fact that Glaisher's earlier result gets to the next power of $p$ compared to Bayat's more recent result above.
\begin{Theorem}(Glaisher, 1900 \cite{GL}) Let $m$ be a positive integer and let $p$ be a prime with $p\geq m+3$. We have,
$$\sum_{k=1}^{p-1}\frac{1}{k^m}\;\equiv\;\begin{cases}\frac{m}{m+1}\,p\,B_{p-1-m}\;\;\qquad\;\;\text{mod}\;p^2\;\;\;\text{if $m$ is even}\\\\
-\,\frac{m(m+1)}{2(m+2)}\,p^2\,B_{p-2-m}\;\;\text{mod}\;p^3\;\;\;\text{if $m$ is odd}
\end{cases}$$
\end{Theorem}

Note, Sun in his proof of $(11)$ and further derivation of the Stirling numbers modulo $p^2$ does not refer to Faulhaber's mathematics, but rather uses his own congruence $(5.1)$, see \cite{SU2}, that is derived from his smart use of Bernoulli's formula.
Be aware of the possible confusion here: what we just referred to as Congruence $(5.1)$ of Sun \cite{SU2} on one hand, and Corollary $(5.1)$ of Sun \cite{SU2} on the other hand are two different things, though both are congruences. From his global formula which reads
\begin{equation}
S_m=p\,B_m+\frac{p^2}{2}\,m\,B_{m-1}+\frac{p^3}{6}\,m(m-1)\,B_{m-2}\;\;\text{mod}\,p^3,
\end{equation}
Sun derives two things using the von Staudt-Clausen theorem. Namely the following two facts: \\\\
\indent $(i)$ $p$ divides $S_m$ as $p$ does not divide the denominator of $B_m$ for $m$ even since $p-1$ does not divide $m$ when $m\leq p-2$.\\\\
\indent $(ii)$ Neither does $p-1$ divide $m-1$ nor $m-2$. Hence, $S_m=p\,B_m\;\text{mod}\,p^2$, independently from the parity of $m$. \\\\
Replacing in $(11)$ yields:
\begin{Result}(Sun 2000, \cite{SU2})
\begin{equation}A_k=\frac{(-1)^{k-1}}{k}\;p\,B_k\;\;\text{mod}\,p^2\qquad\forall\,k=1,\dots,\,p-1\end{equation}
\end{Result}
Unfortunately, Sun's ingenious method for computing the Stirling numbers modulo $p^2$ does not generalize to the modulus $p^3$. Indeed, when working modulo $p^3$, the sum of $(8)$ is no longer congruent to zero, but is rather a convolution of Bernoulli numbers.
We know that $A_{2l+1}=0\;\text{mod}\,p^2$ when $1\leq l\leq \frac{p-3}{2}$ and $p$ divides the $S_j$'s, hence only the terms in even $r$'s contribute to the sum modulo $p^3$, as well as the term corresponding to $r=1$ which is congruent modulo $p^3$ to
$$\frac{p^2}{2}\,B_{k-1}$$ That is, the given sum is congruent modulo $p^3$ to
\begin{equation}\frac{p^2}{2}\,B_{k-1}+\sum_{r=1}^{\lfloor\frac{k-1}{2}\rfloor}A_{2r}\,S_{k-2r}\end{equation}
Moreover, for odd $k$'s, we have when $r\neq \frac{k-1}{2}$ that $p^2$ divides $S_{k-2r}$ as $k-2r$ is then odd and distinct from $1$. Since $p$ also divides $A_{2r}$, the sum to the right hand side of $(26)$ is then congruent modulo $p^3$ to
$$\frac{p^2}{2}\,\frac{B_{k-1}}{k-1},$$
with only the upper index of the sum contributing.\\
On the other hand, for even $k$'s, we have
$$A_{2r}=-p\,\frac{B_{2r}}{2r}\;\text{mod}\,p^2\qquad\text{and}\qquad S_{k-2r}=p\,B_{k-2r}\;\text{mod}\,p^2$$
Thus, modulo $p^3$, equality $(8)$ reads for each $k$ with $1\leq k\leq \frac{p-1}{2}$,
\begin{eqnarray}
(k\leq\frac{p-3}{2})\;A_{2k+1}\negthickspace\negthickspace&=&\negthickspace\!\!\frac{p^2}{2}\frac{2k+1}{2k}\,B_{2k}\qquad\qquad\qquad\qquad\qquad\!\text{mod}p^3\\
&&\notag\\
(k\neq 1)\;A_{2k}\negthickspace\negthickspace&=&\negthickspace\!\!-\frac{1}{2k}\Bigg(p\,B_{2k}-p^2\,\sum_{r=1}^{\lfloor k-\frac{1}{2}\rfloor}\frac{B_{2r}B_{2k-2r}}{2r}\Bigg)\text{mod}p^3\\
&&\notag\\
\text{and}\;A_1\negthickspace\negthickspace&=&\negthickspace\!\!\frac{p(p-1)}{2}\qquad\qquad\qquad\qquad\qquad\qquad\text{mod}\,p^3\\
&&\notag\\
(p\geq 5)\;\text{and}\;A_2\negthickspace\negthickspace&=&\negthickspace\!\!\frac{1}{2}\bigg(-\frac{p}{6}+\frac{3\,p^2}{4}\bigg)\qquad\qquad\qquad\qquad\;\text{mod}\,p^3
\end{eqnarray}
Formula $(27)$ is Result $4$ of Glaisher. Formula $(28)$ is expressed in terms of a convolution involving Bernoulli numbers. Such convolutions have drawn the interest of various mathematicians. For instance, Japanese mathematician Hiroo Miki comes up with an identity in $1978$ which involves both a binomial convolution and an ordinary convolution of divided Bernoulli numbers, see \cite{MI}. Denoting
$$\mathcal{B}_n=\frac{B_n}{n},\;\;\text{and}\;\;H_n=1+\frac{1}{2}+\dots+\frac{1}{n}$$
his identity reads:
$$\forall n>2,\,\sum_{i=2}^{n-2}\mathcal{B}_i\mathcal{B}_{n-i}=\sum_{i=2}^{n-2}\binom{n}{i}\mathcal{B}_i\mathcal{B}_{n-i}+2H_n\mathcal{B}_n$$
Miki shows that both sides of the identity are congruent modulo $p$ for sufficiently large $p$, which implies that they are equal. Another proof of Miki's identity that uses $p$-adic analysis is given by Shiratani and Yokoyama \cite{SH}. Later on, Ira Gessel gives a much simpler proof of Miki's identity based on two different expressions for Stirling numbers of the second kind, \cite{GE}.\\

Unfortunately, we can't resolve $(28)$ by the same method as the one used by Sun in \cite{SU2} in order to compute $(p-1)!$ modulo $p^3$, when he uses the conjugates of the $A_j$'s and of the $S_j$'s. We are thus tempted to apply our own method to the next modulus $p^3$ in order to deal with the Stirling numbers modulo $p^3$. Going to the next modulus $p^3$ happens to provide a much easier way than Sun's for finding $(p-1)!$ modulo $p^3$. It is time to introduce new notations since we work one $p$-power further.
\begin{Notation}
We define $\delta_1(k)$ as the third coefficient in the $p$-adic expansion of $k^{p-1}$, namely,
\begin{equation}
k^{p-1}=1+p\,\delta_0(k)+p^2\,\delta_1(k)\;\;\text{mod}\,p^3
\end{equation}
\end{Notation}
We show the following result.
\begin{Theorem}
\begin{equation}\begin{split}
(p-1)!=-1+p\,\sum_{i=1}^{p-1}\delta_0(i)&+p^2\sum_{i=1}^{p-1}\bigg(\delta_0(i)+\delta_1(i)\bigg)\\
&\\
&-\frac{p^2}{2}\,\Bigg[\bigg(\sum_{i=1}^{p-1}\delta_0(i)\bigg)^2+\sum_{i=1}^{p-1}\delta_0(i)^2\Bigg]\;\;\text{mod}\,p^3
\end{split}\end{equation}
\end{Theorem}
From there, it is easily seen that $(p-1)!$ modulo $p^3$ can be written in terms of sums of powers and powers of sums of powers. Hence, in terms of Bernoulli numbers we get, where we put the terms in random order,
\begin{Corollary}
\begin{equation}(p-1)!=\frac{p}{2}-\frac{3}{2}p^2+(2p+1)pB_{p-1}-\frac{1}{2}p\,B_{2p-2}-\frac{1}{2}\,p^2\,B_{p-1}^2\;\;\text{mod}\,p^3\end{equation}
\end{Corollary}
Sun's formula from \cite{SU2} reads
\begin{equation}(p-1)!=-\frac{p\,B_{p-1}}{p-1}+\frac{p\,B_{2p-2}}{2(p-1)}-\frac{1}{2}\Big(\frac{p\,B_{p-1}}{p-1}\Big)^2\;\;\text{mod}\,p^3\end{equation}
Showing the equivalence between our formula and Sun's formula requests a bit of effort.
First and foremost, we note that
\begin{equation}
(p-1)^{-1}=-p^2-p-1\;\;\text{mod}\,p^3
\end{equation}
Then, by plugging $(35)$ into $(34)$ and using the fact that $p\,B_{2p-2}=-1\;\text{mod}\,p$, namely a consequence of von Staudt-Clausen's theorem, we get,
\begin{equation}
(p-1)!=-2p^2-\frac{p(p+1)}{2}\,B_{2p-2}+p(p+1)\,B_{p-1}-\frac{2p+1}{2}\,p^2B_{p-1}^2\;\;\text{mod}\,p^3
\end{equation}
In order to conclude, we will need to use a consequence of a result of Sun. First we introduce some new notations.
\begin{Notation}
Let $x$ be a $p$-adic integer. We denote by $(x)_k$ the $(k+1)$-th coefficient in the $p$-adic expansion of $x$, that is
$$x=\sum_{j=0}^{\infty}(x)_j\,p^j$$
\end{Notation}
\begin{Fact} By von Staudt--Clausen's theorem, $p\,B_{2p-2}$ and $p\,B_{p-1}$ are both $p$-adic integers of equal residue $-1$ modulo $p$. Their second coefficients of their respective $p$-adic expansions are related by
\begin{equation}\big(p\,B_{2(p-1)}\big)_1=2\;\big(p\,B_{p-1}\big)_1-1\;\;\text{mod}\,p\mb{Z}_p\end{equation}
\end{Fact}
Fact $2$ appears to be a direct consequence of Sun's unpublished result,
\begin{Result} (Sun, \cite{SU3}) Let $k$ be a nonnegative integer. We have
\begin{equation}
p\,B_{k(p-1)}=-(k-1)(p-1)+k\,p\,B_{p-1}\;\;\text{mod}\,p^2
\end{equation}
\end{Result}
Result $6$ of Sun is itself a special case of his Corollary $4.2$ of \cite{SU1}, whose statement we recall here for interest and completeness.
\begin{Result} (Sun, 1997 \cite{SU1})
Let $b$ and $k$ be nonnegative integers. \\
Define $$\delta(n,b,p)=\begin{cases} 1\;\;\text{if $p-1|b$ and $B_n\not\in\mb{Z}_p$}\\
0\;\;\text{otherwise \emph{i.e.} $B_n\in\mb{Z}_p$ or $p-1\not|\,b$}\end{cases}$$
Then, we have
\begin{equation}\begin{split}
\Big(1-p^{k(p-1)+b-1}\Big)&p\,B_{k(p-1)+b}\\&\equiv \sum_{r=0}^{n-1}(-1)^{n-1-r}\binom{k-1-r}{n-1-r}\binom{k}{r}\Big(1-p^{r(p-1)+b-1}\Big)\,p\,B_{r(p-1)+b}\\
&\qquad +(-1)^n\,\delta(n,b,p)\binom{k}{n}\,p^{n-1}\;\;\text{mod}\,p^n\,\mb{Z}_p
\end{split}\end{equation}
\end{Result}
Unpublished Result $6$ was proven prior to more recent Result $7$ and is the latter result applied with $b=0$, $n=2$, $p>3$ and $k\geq 1$.\\ Apply now Result $6$ with $k=2$. Obtain successively
\begin{equation}
p\,B_{2(p-1)}=-(p-1)+2p\,B_{p-1}\;\;\text{mod}\,p^2\,\mb{Z}_p,
\end{equation}
and so \begin{equation}
-1+p\big(p\,B_{2p-2}\big)_1=-p-1+2p\big(p\,B_{p-1}\big)_1\;\;\text{mod}\,p^2\,\mb{Z}_p
\end{equation}
Hence, after simplifying, obtain Fact $2$. From there, concluding is just routine calculation.

\end{itemize}

Our paper is structured as follows. \\

Section $2$ is concerned with proving Theorems $1,2,3$ and Corollaries $1,2$. Some of our proofs, more specifically for Theorem $2$, use Faulhaber's original work from the 17th century, as it is the way we ourselves re-discovered Glaisher's result  from $1900$ (cf Corollary $1$) after we came up with the formula of Theorem $1$ using $p$-adic numbers. Faulhaber didn't get recognition for his pioneering work on sums of powers until much later when mathematicians and computer scientists like Edwards, Gessel, Viennot and Knuth got interested in this branch of mathematics again. They then gave Faulhaber much credit for his long and forgotten groundbreaking work that was done mostly by hand at the time when he lived.
Some $p$-adic properties get derived in our paper from their work. The link of Faulhaber's work with Bernoulli numbers is that there is a uniform formula for the sums of powers of integers in terms of Bernoulli numbers.
This is something that Faulhaber ignored and he had distinguished between odd powers and even powers. The relationship between the coefficients from the old days and the coefficients from the newer age appears in \cite{KN} and is mostly due to \cite{GV}. It is said in \cite{KN} that the book by Faulhaber \cite{JF} is "evidently extremely rare", that "no copies have ever been recorded to exist in America". Edwards found Faulhaber's book at Cambridge University Library, a volume that was once owned by Jacobi. The proof of Theorem $3$ uses Bayat's result (but only modulo $p$) which is a generalization of
Wolstenholme's theorem. A generalization one $p$ power further by Sun for generalized harmonic numbers with even powers is needed in order to derive the result of Corollary $2$. Alternatively, we can use Bayat's result modulo $p^2$ this time, together with Glaisher's way of relating the generalized harmonic numbers to the Stirling numbers modulo $p^2$.

Section $3$ of the paper presents discussions centered around the Agoh-Giuga conjecture and also provides some historical mathematical background. We give a proof for Fact $1$.

In Section $4$, we work modulo $p^3\mb{Z}_p$, still using the same polynomial $g$. First and foremost, we find a formula for $(p-1)!$ mod $p^3$ that is only expressed in terms of sums of powers. Then, we are able to express $(p-1)!$ mod $p^3$ in terms of Bernoulli numbers. Second, we show by using a congruence of Sun type that Sun's formula and our formula are identical.

Finally, Section $5$ summarizes the additions of our work in comparison with Sun and Glaisher's earlier works, evokes the forthcoming plans modulo $p^3$ as announced in the abstract and discusses future research directions and in particular in a quantum setting.

\section{Wilson's theorem modulo $p^2$}

In this part, we first work out the constant coefficient of the polynomial $g$ modulo $p^2\,\mb{Z}_p$ using both the factored form and the expanded form of $g$. This will lead to the formula of Theorem $1$. To that aim, we will first need to investigate the $p$-adic residues of the $t_k$'s which we will denote by $t_k^{(0)}$. Using Hensel's lifting algorithm, we show the following intermediate result.
\newtheorem{Lemma}{Lemma}
\begin{Lemma} Let $k$ be an integer with $1\leq k\leq p-1$. Then, we have:
\begin{equation}p\,t_k^{(0)}=k(1+(p-1)!+p\,\delta_0(k))\qq\end{equation}
\end{Lemma}
\noindent Here, $\delta_0(k)$ is defined as in the statement of Theorem $1$, see $\S\,1$.
We recall from the introduction that $k\in\zsp$ lifts to a unique root $k+pt_k$ of $g$ with $t_k$ a $p$-adic integer. The second coefficient $t_k^{(0)}$ of the root expansion must satisfy to $$g(k+p\,t_k^{(0)})\in\,p^2\mb{Z}_p$$
Hence, we must have
\begin{equation}
k^{p-1}+p(p-1)t_k^{(0)}\,k^{p-2}+(p-1)!\in p^2\mb{Z}_p
\end{equation}
By definition,
$$k^{p-1}=1+p\dk\qquad\text{mod}\,p^2$$
And so,
$$k^{p-2}=k^{-1}+p\dk\,k^{-1}\qq$$
Replacing in $(43)$ yields the equation of Lemma $1$.  \\
Next, by looking at the constant coefficient of $g$ modulo $p^2$, we obtain
\begin{equation}
\sum_{k=1}^{p-1}\,p\tk(p-1)!^k=0\qq
\end{equation}
Then, using $(42)$, it yields
\begin{equation}
(p-1)!\sum_{k=1}^{p-1}\,(1+(p-1)!+p\dk)=0\qq
\end{equation}
We derive
\begin{equation}
(p-1)!\Big((p-1)(1+(p-1)!)+p\sum_{k=1}^{p-1}\dk\Big)=0\qq
\end{equation}
By using Wilson's theorem modulo $p$ inside the bracket, we derive in turn,
\begin{equation}
(p-1)!\Big(-1-(p-1)!+p\sum_{k=1}^{p-1}\dk\Big)=0\qq
\end{equation}
Finally, since $p^2\wedge (p-1)!=1$, we obtain the formula of Theorem $1$. \\The latter formula can be written in terms of a sum of powers of integers as follows,
\begin{equation}
(p-1)!=-p+\sum_{k=1}^{p-1}k^{p-1}\qq
\end{equation}
Assume point $(i)$ of Theorem $2$ holds. Then,
using Faulhaber's formula for sums of even powers as stated in Eq. $(9)$ with $2l=p-1$, we get
\begin{equation}
p\,\sum_{k=1}^{p-1} k^{p-1}=\Big(p-1+\frac{1}{2}\Big)\Big(2\,c_1(l)\frac{p(p-1)}{2}+3\,c_2(l)\frac{p^2(p-1)^2}{4}\Big)\qquad\text{mod}\,p^3
\end{equation}
Further, back in the $17$th century, Faulhaber states that the two trailing coefficients in Eq. $(3)$ will have the form $4\alpha\,a^3-\alpha\,a^2$.
Using this fact and simplifying, we get
\begin{equation}
\sum_{k=1}^{p-1}k^{p-1}=\frac{1}{2}c_1(l)\qquad\text{mod}\,p^2
\end{equation}
Combining $(48)$ and $(50)$ leads to formula $(ii)$ of Theorem $2$. \\
Let's now prove $(i)$. The starting point is Jacobi's formula \cite{JA}. Jacobi was the first to prove Faulhaber's formula in $1834$. Using our former notation for $a$ and letting $u=2a$, Jacobi's formula reads:
\begin{equation}
\sum_{k=1}^{p-1}k^{2l+1}=\frac{1}{p+1}\Big(A_0^{(l+1)}\,u^{l+1}+A_1^{(l+1)}\,u^l+\dots+A_l^{(l+1)}\,u\Big)
\end{equation}
In the most general form of the formula, the denominator $p+1$ should be replaced with $2l+2$. 
Confronting Jacobi's formula with Faulhaber's formula in Eq. $(3)$, we see that
\begin{equation}
A_l^{(l+1)}=0\qquad\text{and}\qquad c_i(l)=\frac{2^{i+1}\,A_{l-i}^{(l+1)}}{p+1}
\end{equation}
It is known that the coefficients $A_k^{(m)}$ obey to some recurrence formulas. Moreover, an explicit formula for these coefficients was first obtained by Gessel and Viennot in \cite{GV} and is provided by Knuth in \cite{KN}. Further, Edwards was first to observe from a recursive formula defining the $A_k^{(m)}$'s and involving binomial coefficients that these numbers can be obtained by inverting a lower triangular matrix, see \cite{Ed}. From there, Gessel and Viennot expressed the coefficients in terms of a $k\times k$ determinant, namely
\begin{equation*}
A_{k}^{(m)}=\frac{1}{(1-m)\dots (k-m)}\left|\begin{array}{ccccc}\bn{m-k+1}{3}&\bn{m-k+1}{1}&0&\dots&0\\
\bn{m-k+2}{5}&\bn{m-k+2}{3}&\bn{m-k+2}{1}&\dots&0\\
\vdots&\vdots&\vdots&&\\
\bn{m-1}{2k-1}&\bn{m-1}{2k-3}&\bn{m-1}{2k-5}&\dots&\bn{m-1}{1}\\
\bn{m}{2k+1}&\bn{m}{2k-1}&\bn{m}{2k-3}&\dots&\bn{m}{3}
\end{array}\right|
\end{equation*}
Of interest to us here, we note that this determinant is an integer. In fact, this determinant has a neat combinatorial interpretation that is due to Gessel and Viennot. A result in \cite{GV} states that this determinant counts the number of ways to put positive integers into a $k$-rowed tripled staircase with numbers strictly increasing from left to right and from top to bottom and imposing also that an entry in row $j$ is at most $m-k+j$.
\begin{center}\epsfig{file=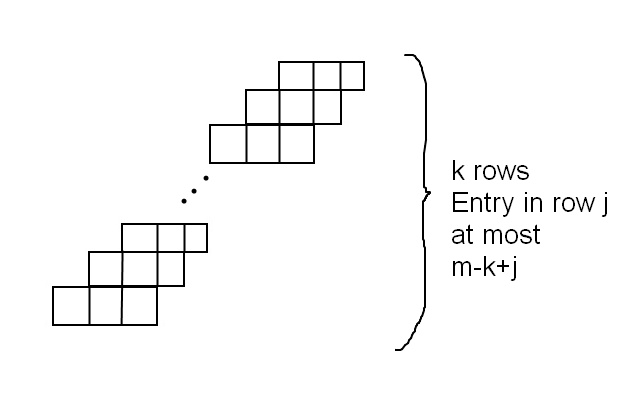, height=5cm}\end{center}
Our goal is to show that $c_i(l)\in\mb{Z}_p$. For that, we use formula $(52)$. From the discussion above, we know that
$$A_{l-i}^{(l+1)}=\frac{1}{l(l-1)(l-2)\dots (l-(l-i-1))}\times \text{an integer}$$
Now $l=\frac{p-1}{2}<p$, hence $p$ does not divide the denominator of $A_{l-i}^{(l+1)}$. Thus, $A_{l-i}^{(l+1)}$ is a $p$-adic integer and so is $c_i(l)$ for each $i$. This closes the proof of Theorem $2$. We are now in a position to retrieve Glaisher's result from $1900$. By \cite{KN},
$$A_{m-2}^{(m)}=\bn{2m}{2}\,B_{2m-2},\;\;\forall\,m\geq 2,$$
hence \begin{equation}c_1(l)=\frac{4A_{l-1}^{(l+1)}}{p+1}=\frac{4p(p+1)}{2(p+1)}\,B_{p-1}=2p\,B_{p-1}\end{equation}
By $(48)$, $(50)$ and $(53)$, we conclude that
$$(p-1)!=p\,B_{p-1}-p\qq$$
The method just exposed allows to bypass the use of von Staudt-Clausen's theorem. Another way namely of getting Glaisher's formula from Theorem $1$ would be to use Bernoulli's formula together with von Staudt-Clausen's theorem as follows. \\
A trick consists of summing the $(p-1)$th powers of integers up to $p$ instead of $p-1$ since we work modulo small powers of $p$ and the last term of the sum will not contribute anyway. This way, Bernoulli's formula can be expressed in terms of powers of $p$ instead of powers of $(p-1)$.
We get
\begin{equation}
p\,\sum_{k=1}^{p}k^{p-1}=\sum_{j=0}^{p-1}\bn{p}{j}B_j\,p^{p-j}
\end{equation}
By $(48)$, we need to know $p\,\sum_{k=1}^p\,k^{p-1}$ modulo $p^3$. This is why we need to know the $p$-divisibility properties of the denominators of the Bernoulli numbers. By von Staudt-Clausen's theorem, amongst the Bernoulli numbers present in the sum of Eq. $(54)$, $p$ divides only the denominator of $B_{p-1}$. We thus deduce
\begin{equation}
p\,\sum_{k=1}^{p-1}k^{p-1}=p^2\,B_{p-1}\qquad\text{mod}\,p^3\mb{Z}_p
\end{equation}
And so,
\begin{equation}
\sum_{k=1}^{p-1}k^{p-1}=p\,B_{p-1}\qquad\text{mod}\,p^2\mb{Z}_p
\end{equation}
Glaisher's formula follows from Eqs. $(48)$ and $(56)$.\\
We note from Glaisher's formula an improvement for calculating $(p-1)!$ modulo $p^2$. For instance, the one millionth Bernoulli number has
$4767554$ digits over $24$ digits and the $1.5$ millionth Bernoulli number has $7415484$ digits over $55$ digits.
A work by Derby dating from $2015$ has provided an efficient way for computing $p\,B_{p-1}$. We must introduce new notations.
\begin{Notation}
Denote by $\widehat{PT}(n)$ the lower triangular matrix obtained by entering the coefficients of the Pascal triangle with the ending $1$ coefficient omitted.
\end{Notation}
Write \begin{eqnarray*}\sum_{k=1}^{p-1}k^p&=&\frac{1}{p+1}\sum_{j=0}^{p}\binom{p+1}{j}B_j\,(p-1)^{p+1-j}\\&=&d_1(p-1)+d_2(p-1)^2+\dots+d_{p+1}(p-1)^{p+1},\end{eqnarray*}
according to Bernoulli's formula.
Derby's result provides a way of computing the coefficients $d_i$'s by simply inverting a matrix involving the Pascal triangle.
\begin{Fact} (Due to Derby \cite{NI})
We have $$(1\;p\;\dots\;p\;1)\widehat{PT}(p+1)^{-1}=(d_1,\,d_2,\,\dots,\,d_{p+1}),$$
with the row to the left hand side denoting the $p$th row of the Pascal triangle.
\end{Fact}
From there, since the second coefficient is
$$\frac{1}{p+1}\binom{p+1}{p-1}\,B_{p-1}=\frac{p}{2}\,B_{p-1},$$
it provides an efficient way of computing $(p-1)!$ modulo $p^2$. \\For interest and completeness here, we note that similar works as Derby's and additional results involving the sums of powers and the Pascal triangle were achieved by Italian mathematician Giorgio Pietrocola and published electronically in $2017$ under the title "On polynomials for the calculation of sums of powers of successive integers and Bernoulli numbers deduced from the Pascal's triangle".

As announced in the introduction, our method also provides interesting congruences modulo $p^2$ on the Stirling numbers. Namely,
by looking at the coefficients of powers of $X$ modulo $p^2$, we can prove Theorem $3$. First, we study the even powers. We obtain the following congruence, where we used Notation $1$:
\begin{multline}
\left[\begin{array}{l}\;\;\;\,p\\2n+1\end{array}\right]+\\\sum_{m_0}\frac{p\,t_{m_0}^{(0)}(p-1)!^{m_0}}{(2n)!}\sum_{m_1\neq m_0}\frac{1}{m_1}\sum_{m_2\neq m_0,m_1}\frac{1}{m_2}\dots\sum_{m_{2n}\neq m_0,m_1,\dots,m_{2n-1}}\frac{1}{m_{2n}}\\=0\,\text{mod}\,p^2
\end{multline}
As $2n\leq p-3$, it suffices to resolve the series of sums on $m_1,\dots,m_{2n}$ modulo $p$. We resolve the sums from right to left, summing over all the terms first and subtracting the illicit terms next. By using Result $1$ of Bayat, we get
$$\sti=-(p-1)!\,\sum_{m_0=1}^{p-1}\frac{p\,t_{m_0}^{(0)}}{m_0^{2n+1}}\;\text{mod}\,p^2$$
Then, it follows from using Lemma $1$,
$$\sti=-(p-1)!\,\sum_{m_0=1}^{p-1}\frac{1+(p-1)!+p\delta_0(m_0)}{m_0^{2n}}\;\text{mod}\,p^2$$
After applying Wilson's theorem modulo $p$ to the factorial to the left, we get
$$\sti=(1+(p-1)!)\sum_{m_0=1}^{p-1}\frac{1}{m_0^{2n}}+\sum_{m_0=1}^{p-1}\frac{p\delta_0(m_0)}{m_0^{2n}}\;\text{mod}\,p^2$$
The firs term is a product of two terms that are congruent to zero modulo $p$ by Wilson's theorem and Bayat's result respectively.
As for the second term, it can be rewritten as a difference after applying the equality
$$p\delta_0(m)=m^{p-1}-1\;\text{mod}\,p^2$$
We obtain
\begin{equation}\sti=\sum_{m=1}^{p-1}m^{p-1-2n}-\sum_{m=1}^{p-1}\frac{1}{m^{2n}}\;\;\text{mod}\,p^2\end{equation}
We get the formula of Theorem $3$ for odd indices $k$. \\
By looking at the coefficient of $X^{2n+1}$ in both factored and expanded forms, we get by the same methodology as before,
\begin{eqnarray*}\stp&=&(p-1)!\,\sum_{l}\frac{p\,t_l^{(0)}}{l^{2n+2}}\;\text{mod}\,p^2,\\
&=&-\sum_{l=1}^{p-1}\frac{p\,t_l^{(0)}}{l^{2n+2}}\;\text{mod}\,p^2\\
&=&-\sum_{l=1}^{p-1}\frac{(1+(p-1)!+p\,\delta_0(l))}{l^{2n+1}}\;\text{mod}\,p^2
\end{eqnarray*}
From there, we derive in a similar way as before using Bayat's result that $\sum_{l=1}^{p-1}\frac{1}{l^{2n+1}}=0\,\text{mod}\,p$,
\begin{equation}\stp=-\Bigg(\sum_{l=1}^{p-1}l^{p-2n-2}-\sum_{l=1}^{p-1}\frac{1}{l^{2n+1}}\Bigg)\;\text{mod}\,p^2\end{equation}
Up to a minus sign, this is like in the statement of Theorem $3$, for even $k$ this time. In fact, as we will see below, the right hand side of Eq. $(59)$ is zero. Then, the congruence of Theorem $3$ holds anyhow and it thus closes the proof of Theorem $3$. \\
By using formula $(24)$ and von Staudt-Clausen's theorem, we have
$$S_{p-k}=p\,B_{p-k}\;\text{mod}\,p^2$$
Moreover, by the top congruence of Sun's Corollary $5.1$ of \cite{SU2}, we have for every integer $k$ with $2\leq k\leq p-1$,
$$H_{p-1,k-1}=\frac{k-1}{k}\,p\,B_{p-k}\;\text{mod}\,p^2$$
Further, when $k$ is even, $p-k$ is odd, hence $B_{p-k}=0$ except when $p-k=1$, that is when $k=p-1$. That case gets treated as a direct calculation. We thus obtain the result of Corollary $2$. \\
We now present another method for deriving Corollary $2$ from Theorem $3$, one which does not use Sun's congruences for generalized harmonic numbers modulo $p^2$, nor Sun's formula for the sums of powers modulo $p^2$. First, we show that $p^2$ divides the Stirling numbers with even indices, to the exception of index $(p-1)$. Going back to congruence $(59)$, we show that the two sums that are involved are both congruent to zero modulo $p^2$. The fact that this holds for the first sum follows directly from applying Faulhaber's formula for odd powers, joint with our point $(i)$ of Theorem $2$ which still holds when $p$ is being replaced with any $2l+1$ such that $2l+1\leq p-2$. Indeed, we have
$$\forall 1\leq i\leq l-1,\;c_i(l)=\frac{2^{i+1}\,A_{l-i}^{(l+1)}}{2l+2}\qquad\text{and}\qquad A_0^{(l+1)}=1$$
From the work of Gessel and Viennot, when $1\leq i\leq l-1$, the denominator of $A_{l-i}^{(l+1)}$ is, up to its sign, a product of integers ranging from $i+1$ to $l$. And so, all the Faulhaber coefficients are $p$-adic integers when the odd power is less than or equal to $p-2$ (note this is no longer necessarily true in general, but the fact still holds when the odd power is $p$, as we proved earlier). As for the second sum, since $2n+1\leq p-4$, an application of Hardy and Wright page $103$ of \cite{HW} provides the desired congruence. Once we know this extra divisibility property on the Stirling numbers, a smart application of Newton's formulas by Glaisher allows to relate the generalized harmonic numbers to the Stirling numbers either modulo $p^2$ or modulo $p^3$. This is Theorem $4$ of the introduction. In particular, from Glaisher's results, we have for odd $k$, using our own notations, 
$$H_{p-1,k-1}=(k-1)\left[\begin{array}{l}p\\k\end{array}\right]\;\;\text{mod}\,p^2$$
Then, replacing in $(58)$ yields the congruence for odd $k$,
$$k\,\left[\begin{array}{l}p\\k\end{array}\right]=S_{p-1,p-k}\;\;\text{mod}\,p^2,$$
thus leading ultimately to Corollary $2$ by a different proof.

\section{Around the Agoh-Giuga conjecture}
We have seen earlier in $\S\,1$ that two types of composite numbers play a central role around the Agoh-Giuga conjecture, namely the Carmichael numbers and the Giuga numbers. For the definition of these numbers, we refer the reader to the introduction. By a theorem of Giuga, the Agoh-Giuga conjecture states that there does not exist any composite number which is both a Carmichael number and a Giuga number. \\It follows from the definition that a Carmichael number is odd and has at least three prime factors. Indeed, notice that
if $n=p\,u$ with $p$ prime, then
$$p-1|n-1\;\Leftrightarrow\; p-1|u-1$$
as $$(n-1)-(u-1)=u(p-1)$$
In other words, $$p-1|n-1\;\Leftrightarrow\;p-1|\frac{n}{p}-1\qquad (E)$$
Consequently, if $2$ divides a Carmichael number and $p$ is another prime divisor of that Carmichael number, then the even number $p-1$ divides an odd number, which is impossible. Thus, a Carmichael number must be odd. \\
Next, suppose for a contradiction that $n$ is Carmichael with $n=pq$ with $p$ and $q$ primes and $p>q$. Since $p-1>q-1$, we have $p-1\not|\;q-1$. But by $(E)$, we have $p-1|\frac{n}{p}-1=q-1$, leading to a contradiction.\\
It also follows from the definition that a Giuga number is square-free, as otherwise, both $\frac{n}{p}$ and $\frac{n}{p}-1$ would be divisible by $p$ when $p$ is a prime dividing $n$. This is impossible.
Thus, we see with the latter equivalence $(E)$ that the Agoh-Giuga conjecture reformulates as:
"There does not exist any composite $n$ such that for all prime divisor $p$ of $n$, the integer $\frac{n}{p}-1$ is divisible by both $p$ and $p-1$".\\
The first seven Carmichael numbers were discovered by Czech mathematician V$\overset{\vee}{a}$clav $\overset{\vee}{S}$imerka in $1885$. They are: \begin{eqnarray*}
561&=& 3.11.17\\
1105&=& 5.13.17\\
1729&=&7.13.19\\
2465&=&5.17.29\\
2821&=&7.13.31\\
6601&=&7.23.41\\
8911&=&7.19.67
\end{eqnarray*}
In $1994$, Alford, Granville and Pomerance show that there exist infinitely many Carmichael numbers. Their proof appears in \cite{AL}.
It is unknown whether there are finitely or infinitely many Giuga numbers. The following data on Giuga numbers are taken from the quite complete and interesting exposition of \cite{BO}. The smallest Giuga number are:
\begin{eqnarray*}
30:&&\frac{1}{2}+\frac{1}{3}+\frac{1}{5}-\frac{1}{30}=1\\
858:&&\frac{1}{2}+\frac{1}{3}+\frac{1}{11}+\frac{1}{13}-\frac{1}{858}=1\\
1722:&&\frac{1}{2}+\frac{1}{3}+\frac{1}{7}+\frac{1}{41}-\frac{1}{1722}=1
\end{eqnarray*}
The smallest odd Giuga number has at least $9$ prime factors since with a smaller number of prime factors, the sum $\frac{1}{p_1}+\dots+\frac{1}{p_m}-\frac{1}{n}$ is smaller than $1$. Whether an odd Giuga number does actually exist or not remains an open question. If we could show that there does not exist any such number, then the Agoh-Giuga conjecture would be proven by Giuga's theorem since a Carmichael number must be odd.\\
Giuga used the property that if $p$ is a prime factor of a Carmichael number $n$, then for no $k$ is $kp+1$ a prime factor of $n$ (as otherwise, $kp+1-1=kp|n-1$ and $p$ divides $n$, a contradiction), to prove computationally that any counterexample to the conjecture would have at least $1000$ digits. Bedocchi in \cite{BED} furthered the number of digits to $1700$ by using the same method. The authos of \cite{BO} improved on this method by reducing the number of cases to be looked at and have shown computationally that any counterexample has no less than $13800$ digits. They also introduce what they call "Giuga sequences" by no longer considering primes in Definition $1$, but rather considering increasing sequences of integers. Like there is an interpretation for a Giuga number in terms of divisibility property, there is a similar approach for Giuga sequences. Inspired by that fact, they also define Carmichael sequences and show that Giuga's conjecture would be proved if one were to show that no Giuga sequence can be a Carmichael sequence. They list some open questions concerning Giuga's sequences and Giuga's conjecture. \\

We now present some of our own thoughts on this interesting topic.
If the Agoh-Giuga conjecture did not hold, there would exist a composite $n$ such that $n\,B_{n-1}=-1\;\text{mod}\,n$. In particular, $n$ is odd. Also, $n\,B_{n-1}\neq 0\,\text{mod}\,n$. As we have seen in the introduction of the current paper, this implies by von Staudt-Clausen's theorem that $n$ has a prime divisor $p$ such that $n$ is Carmichael at $p$. In fact, we show that the converse holds if we also impose that $n$ is square-free.
\newtheorem{Proposition}{Proposition}
\begin{Proposition}
Let $n$ be composite and square-free. Then,
$$n\,B_{n-1}\neq 0\;\text{mod}\,n\;\Leftrightarrow\exists\,p|n\;\text{such that}\;p-1|n-1$$
\end{Proposition}
In order to prove Proposition $1$, we first state the following lemma.
\begin{Lemma}
$$p\,B_{n-1}=\begin{cases} 0&\text{mod}\,p\qquad\text{if}\;p-1\not|\,n-1\\ -1&\text{mod}\,p\qquad\text{if}\;\;p-1\;|n-1\end{cases}$$
\end{Lemma}
Lemma $2$ follows from the fact that $p\,B_{k(p-1)}=-1\;\text{mod}\,p$. This is simply a consequence of von Staudt-Clausen's theorem.
Let $n$ be composite and square-free such that $p$ is a prime number dividing $n$. Then, we may write
$$n=p\prod_{\begin{array}{l}q_i\in\mathcal{P}\\q_i\neq p\\q_i|n\end{array}}q_i$$
It follows from Lemma $2$ that
$$n\,B_{n-1}=\begin{cases} 0&\text{mod}\,p\qquad\text{if}\;p-1\not|\,n-1
\\
-\underset{\begin{array}{l}q_i\in\mathcal{P}\\q_i\neq p\\q_i|n\end{array}}{\prod}q_i&\text{mod}\,p\qquad\text{if}\;\;p-1\;|n-1\end{cases}$$
Thus, if $p-1|n-1$, then $n\,B_{n-1}\neq 0\;\text{mod}\,p$ and consequently $n\,B_{n-1}\neq 0\;\text{mod}\,n$.
We deduce a corollary.
\begin{Corollary} Let $n$ be composite and square-free. Suppose $n$ is Carmichael at at least one of its prime divisors. Then,
at least one of these primes $p$ at which $n$ is Carmichael must satisfy to $p\not|\;\;N(B_{n-1})$.
\end{Corollary}
The proof of Corollary $7$ is as follows. By Proposition $1$, we know that $n\,B_{n-1}\neq 0\;\text{mod}\,n$. If for all the prime factors of the denominator of $B_{n-1}$, these primes also divide the numerator of $B_{n-1}$, then $B_{n-1}$ is an integer. Then $n\,B_{n-1}=0\;\text{mod}\,n$, which constitutes a contradiction. \\

Before moving further, it is a good time to introduce irregular primes and their role in the numerators of the divided Bernoulli numbers.
\begin{Definition} Let $p$ be a prime number. \\We say that
$p$ is a regular prime if $p$ does not divide any of $B_2$, $B_4$, $\dots$, $B_{p-3}$. \\The prime $p$ is said to be irregular otherwise.
\end{Definition}
Regular primes were introduced by K\"ummer when he proved that Fermat's last theorem is true when the exponent is a regular prime. It is unknown whether there exist infinitely many regular primes. On the other hand, Jensen proved in $1915$ that the number of irregular primes is infinite \cite{JEN}. He more precisely shows that there are infinitely many irregular primes that are congruent to $3$ mod $4$. There are only three irregular primes below $100$: $37$, $59$ and $67$. Part of Tanner's undergraduate thesis at Harvard under John Tate searches for irregular primes $p$ with $125\,000<p<150\,000$. The results appear in \cite{TW}. There are only four such primes which each divide four Bernoulli numbers, that is there are four associated irregular pairs for each of these four irregular primes. We say that the index of irregularity of these primes is $4$. Tanner and Tate proved Fermat's last theorem for these irregular primes. Previously in \cite{WA}, Wagstaff had determined the irregular primes up to $125\,000$. Other known facts on irregular primes include the fact that the only irregular pairs $(p,p-3)$ for $p<10^9$ are obtained for $p=16\,843$ (result of Selfridge and Pollack \cite{SE}) and for $p=21\,24\,679$ (result of Buhler, Crandall, Ernvall and Mets\"anky\"a \cite{BU}). Such a prime is called a Wolstenholme prime as for these primes, $p$ divides the Wolstenholme quotient, namely
\begin{eqnarray*}\binom{2p-1}{p-1}=1\;\text{mod}\,p^4\Leftrightarrow p\,\,|W_p=\frac{\binom{2p-1}{p-1}-1}{p^3}&\Leftrightarrow& H_{p-1,1}=0\;\text{mod}\,p^3\\&\Leftrightarrow& H_{p-1,2}=0\;\text{mod}\,p^2\end{eqnarray*}
In the general case, Glaisher had shown that $W_p=-\frac{2}{3}\,B_{p-3}\,\text{mod}\,p$.

Like mentioned at the beginning of the paragraph, the main interest of distinguishing between regular primes and irregular primes when dealing with Bernoulli numbers relies on the following fact.
\begin{Fact}
$$Numer\Bigg(\frac{B_{n}}{n}\Bigg)=\begin{cases} 1 \;\;\text{when $n=2,4,6,8,10,14$}\\
\text{A product of powers of irregular primes otherwise}
\end{cases}$$
\end{Fact}

In what follows, we will make use of K\"ummer's congruences which we recall below.
\begin{Theorem}(K\"ummer's congruences, $1850$ \cite{KU})\\
Let $p$ be an odd prime and let $b>0$ be an even integer such that $p-1\not|\,b$. Then, we have for all nonnegative integer $k$,
$$\frac{B_{k(p-1)+b}}{k(p-1)+b}\equiv\frac{B_b}{b}\;\;\text{mod}\,p$$
\end{Theorem}
Note, in \cite{SU2}, Sun generalizes these congruences to the moduli $p^2$ and $p^3$. \\
Carlitz gave a proof that there exist infinitely many irregular primes, using K\"ummer's congruences \cite{CAR2}. We outline his proof below. \\
Let $\lbrace p_1,\dots,p_s\rbrace$ be the set of irregular primes and set $n=k(p_1-1)\dots(p_s-1),$ with $k$ chosen even and so that
$$v_p\Big(\frac{B_n}{n}\Big)>0$$ for some prime $p$ (recall that $|\frac{B_{2m}}{2m}|\underset{m\to +\infty}{\longrightarrow}+\infty$).\\
If $p-1|n$, then $n=(p-1)l$ for some integer $l$. By von Staudt-Clausen's theorem, $pB_{l(p-1)}=-1\;\text{mod}\,p$. But since
$$\frac{B_{l(p-1)}}{l(p-1)}=p\,u,$$ for some $u\in\mb{Z}_p$, we see that $p\,B_{l(p-1)}=p^2\,u\,l(p-1)=0\,\text{mod}\,p$, a contradiction.
Hence $p-1\not|\,n$ and in particular $p$ is none of the $p_i$'s. We conclude by showing that $p$ is irregular.
Write $n=(p-1)s+t$ with $0<t<p-1$. We notice that $t$ is even as $n$ is even. Moreoever, $p-1\not|\,t$. Then, we are under the conditions of application of K\"ummer's congruences. It yields:
$$\frac{B_n}{n}\equiv\frac{B_t}{t}\;\;\text{mod}\,p$$
This forces that $p|B_t$ with $t\leq p-3$. Thus, $p$ is irregular. \\

In light of K\"ummer's congruences, we are now prepared for proving Fact $1$ of the introduction. Fact $1$ asserts that if $n$ is an odd composite Giuga number, then for a divisor $p$ of $n$, either $n$ is Carmichael at $p$ or $p$ divides $B_{\frac{n}{p}-1}$. \\We prove this fact below.

Let $n=m\,p$ with $n$ Giuga and odd. In particular $n$ is square-free (from being Giuga). Since $n$ is Giuga, we have by definition $p|m-1$. Since $n$ is odd, $m$ must be also odd, and so $m-1$ is even.
Suppose $n$ is not Carmichael at $p$, that is $p-1\,\not|\;m-1$ by $(E)$. Then, an application of Theorem $7$ yields for all nonnegative integer $k$ :
$$\frac{B_{m-1}}{m-1}\equiv\frac{B_{k(p-1)+m-1}}{k(p-1)+m-1}\;\;\text{mod}\,p$$
Pick $k$ such that $p\not|k$. Then $v_p(m-1+k(p-1))=0$, where $v_p$ denotes the $p$-adic valuation.
Also, since $p-1$ does not divide $m-1+k(p-1)$, we know by von Staudt-Clausen's theorem that $p$ does not divide $D(B_{m-1+k(p-1)})$. Then the right hand side of the latter congruence is a $p$-adic integer. Then, the congruence itself imposes that the left hand side of the congruence is also a $p$-adic integer. However, $p$ divides $m-1$. It then forces $p|B_{m-1}$, just like expected.\\

\noindent We note that Fact $1$ appears as a straightforward consequence of Adams'theorem. We recall Adams'statement below.
\begin{Theorem} (Adams'theorem \cite{AD} $1878$)
Let $p$ be an odd prime with $p>3$. \\
If $p^l|n$ and $p-1\not|\,n$, then $p^l|B_n$
\end{Theorem}
Note, Adams'theorem does not apply when $p=3$. But when $p=3$, since $n$ is odd, $\frac{n}{3}-1$ is even and so $2$ divides $n-1$. Thus, $n$ is Carmichael at $3$ and Fact $1$ holds. \\
There exists a more recent version of Adams theorem stated by Thangadurai as a conjecture. By $p^a||n$, he means that $p^a|n$, but $p^{a+1}\not|,n$.
\begin{Conjecture}(Thangadurai \cite{TA}, $2004$)
Let $p>3$ be a prime such that $p^l||n$ for some even positive integer $n$ and $p-1\not|,n$. \\
Then, $p^{\beta}||N(B_n)$ implies that $\beta\leq l+1$
\end{Conjecture}
Thangadurai shows that his conjecture holds for regular primes and for irregular primes less than $12$ millions.
While Adams does not provide any proof for his theorem in his own paper, there exist a few proofs, see \cite{JO} and \cite{IR}. None of these two proofs is based on K\"ummer's congruences. The proof of \cite{IR} uses Voronoi's congruences.

\section{Wilson's theorem modulo $p^3$}

Like we did in the modulus $p^2$ case, we must lift the root residues of $g$ one $p$ power further. First, we introduce a new notation.
\begin{Notation}
Define $t_k^{(1)}$ as the $p$-adic residue such that
$$x_k=k+p\,t_k^{(0)}+p^2\,t_k^{(1)}\qquad\text{mod}\,p^3$$
\end{Notation}
\noindent The following lemma will be useful to the proof of Theorem $6$. We used Notation $2$ from $\S\,1$.
\begin{Lemma}
$$t_k^{(1)}=k\Big(\delta_0(k)+\delta_1(k)+\Bigg(\sum_{i=1}^{p-1}\delta_0(i)\Bigg)^2+(1+\delta_0(k))\,\sum_{i=1}^{p-1}\delta_0(i)\Big)\;\;\text{mod}\,p$$
\end{Lemma}
\noindent Hensel's lifting algorithm is so that
$$g(k+p\,t_k^{(0)}+p^2\,t_k^{(1)})\in\,p^3\mb{Z}_p,$$
that is
$$(k+p\,t_k^{(0)}+p^2\,t_k^{(1)})^{p-1}+(p-1)!\in\,p^3\mb{Z}_p$$
Expanding yields:
\begin{multline*}k^{p-1}+p(p-1)\tk k^{p-2}+\frac{(p-1)(p-2)}{2}p^2\big(\tk\big)^2 k^{p-3}\\+p^2(p-1)\tku (k+p\,\tk)^{p-2}+(p-1)!\in p^3\,\mb{Z}_p\end{multline*}
Another round of simplifications modulo $p^3$ lead to
$$k^{p-1}+p\,\tk(p-1)k^{p-2}-p^2\tku k^{p-2}+p^2\big(\tk\big)^2\,k^{p-3}+(p-1)!\in p^3\,\mb{Z}_p$$
Notice that since $0\leq\tk\leq p-1$, we have $\tk\;\text{mod}\,p^2=\tk\;\text{mod}\,p$. It follows that
$p\tk\;\text{mod}\,p^3=p\,\tk\;\text{mod}\,p^2$. Thus, we have by Lemma $1$,
$$p\tk=k(1+(p-1)!+p\,\dk)\;\;\text{mod}\,p^3$$
We get
\begin{multline*}p^2\,t_k^{(1)}\,k^{-1}=k^{p-1}\Big(1+(p-1)(1+(p-1)!+p\dk)+(1+(p-1)!+p\dk)^2\Big)\\
+(p-1)!\;\;\text{mod}\,p^3\end{multline*}
After replacing $$k^{p-1}=1+p\dk+p^2\delta_1(k)\;\;\text{mod}\,p^3,$$ using Theorem $1$ when appropriate and simplifying modulo $p^3$, we obtain the expression of the lemma. Theorem $6$ of $\S\,1$ is then derived by looking at the constant coefficient of $g$ modulo $p^3$ in both factored and expanded forms. We have, using our Notation $1$
\begin{multline*}(p-1)!=(p-1)!+\sum_{k=1}^{p-1}p\,\tk\,(p-1)!^k+\sum_{k=1}^{p-1}p^2\,t_k^{(1)}(p-1)!^k\\
+\sum_{i\neq j}pt_i^{(0)}pt_j^{(0)}(p-1)!^{i,j}\;\;\text{mod}\,p^3
\end{multline*}
Then, factoring $(p-1)!$ and using Lemma $1$ and Lemma $3$, we get
\begin{multline*}
(p-1)!\Bigg(\sum_{k=1}^{p-1}(1+(p-1)!+p\dk)+p^2\s(\dk+\dku)\\-p^2\bigg(\si\di\bigg)^2+p^2\bigg(\si\di\bigg)\s(1+\dk)\\
+\sum_{i\neq j}\big(1+(p-1)!+p\di\big)\big(1+(p-1)!+p\djz\big)\Bigg)=0\;\;\text{mod}\,p^3
\end{multline*}
We derive,
\begin{multline*}
(p-1)!\Bigg(p^2\si\di-1-(p-1)!+p\si\di+p^2\si(\diu+\di)\\
-p^2\bigg(\si\di\bigg)^2-p^2\si\di+p^2\bigg(\si\di\bigg)^2\\
+\frac{1}{2}\si\sum_{j\neq i}\big(1+(p-1)!+p\di\big)\big(1+(p-1)!+p\djz\big)\Bigg)=0\mdt
\end{multline*}
In the left hand side above, some terms simplify. Denote by $S$ the double sum.
We evaluate it as follows.
\begin{multline*}
S=\frac{1}{2}\Bigg(\si(1+(p-1)!+p\di)\sum_{j=1}^{p-1}(1+(p-1)!+p\djz)\\-\si\bigg(1+(p-1)!+p\di\bigg)^2\Bigg)\mdt
\end{multline*}
A quick inspection shows that the first term of the difference above is zero modulo $p^3$. Thus, modulo $p^3$, the double sum $S$ reduces to
$$S=-\frac{p^2}{2}\Bigg(\bigg(\si\di\bigg)^2+\si\di^2\Bigg)\mdt$$
By gathering the different parts, we obtain the formula of Theorem $6$. It remains to express $(p-1)!$ modulo $p^3$ in terms of sums of powers and powers of sums of powers. From there, we will derive the formula of Corollary $6$, that is Wilson's theorem modulo $p^3$. We prove the following lemma.
\begin{Lemma}
\begin{multline*}
(p-1)!=-1-\frac{1}{2}\Bigg(\bigg(\si i^{p-1}\bigg)^2+\si(i^{p-1})^2\Bigg)+(2p+1)\si i^{p-1}\\
-(p-1)(\frac{3}{2}p+1)\mdt
\end{multline*}
\end{Lemma}
\noindent Start from the formula of Theorem $6$. First and foremost, we group the terms so as to use the expansion
$$i^{p-1}=1+p\di+p^2\diu\mdt$$
We get,
\begin{multline}
(p-1)!=-1+\si (p\di+p^2\diu)+p^2\si\di\\
-\frac{1}{2}\si(i^{p-1}-1-p^2\diu)^2-\frac{1}{2}\Bigg(\si i^{p-1}-(p-1)-p^2\si\diu\Bigg)^2
\end{multline}
Denote the last two sums of the latter expression respectively by $S_3$ and $S_4$. We proceed to the evaluation of these two sums.
We have by expanding the square of $S_3$,
$$S_3=-\frac{1}{2}\si i^{2p-2}-\frac{1}{2}(p-1)+\si i^{p-1}+p^2\si \diu\bigg(i^{p-1}-1\bigg)$$
Then, using the expression for $p^2\diu$ modulo $p^3$ and replacing, it comes
$$S_3=-\frac{1}{2}\si i^{2p-2}-\frac{1}{2}(p-1)+\si i^{p-1}+\si (i^{p-1}-1-p\di)(i^{p-1}-1)$$
By expanding the factor of the last sum and regrouping the terms, it follows that
$$S_3=\frac{1}{2}\si i^{2p-2}+\frac{1}{2}(p-1)-\si i^{p-1}+p\si\di-p\si\di i^{p-1}$$
Now, notice that
$$p\di\,p\di=p\di\,(i^{p-1}-1)\mdt$$
Hence, the combination of the last two terms of the latter expression giving $S_3$ is nothing else than $2\,S_3$. It follows that
\begin{equation}
S_3=-\frac{1}{2}\,\si i^{2p-2}+\si i^{p-1}-\frac{1}{2}(p-1)\mdt
\end{equation}
We now deal with $S_4$. Firs, we expand the square. It yields:
\begin{multline*}S_4=-\frac{1}{2}\bigg(\si i^{p-1}\bigg)^2-\frac{1}{2}(p-1)^2+(p-1)\si i^{p-1}\\+\bigg(p^2\si \diu\bigg)\bigg(\si i^{p-1}-(p-1)\bigg)\mdt\end{multline*}
Next, by replacing $p^2\diu$, it comes:
\begin{multline*}S_4=-\frac{1}{2}\bigg(\si i^{p-1}\bigg)^2-\frac{1}{2}(p-1)^2+(p-1)\si i^{p-1}\\+\Bigg(\si\big(i^{p-1}-1-p\di\big)\Bigg)\Bigg(\si i^{p-1}-(p-1)\Bigg)\mdt\end{multline*}
Expanding the last term and simplifying leads in turn to
$$S_4=\frac{1}{2}\bigg(\si i^{p-1}\bigg)^2+\frac{1}{2}(p-1)^2-(p-1)\si i^{p-1}+2\,S_4\mdt,$$
from which we finally derive
\begin{equation}
S_4=-\frac{1}{2}\bigg(\si i^{p-1}\bigg)^2-\frac{1}{2}(p-1)^2+(p-1)\si i^{p-1}\mdt
\end{equation}
By plugging the respective expressions $(61)$ and $(62)$ for $S_3$ and $S_4$ into $(60)$ and simplifying, we obtain the formula of Lemma $4$.\\
From there, concluding to Wilson's theorem modulo $p^3$ as in our Corollary $6$ is simply a matter of applying von Staudt-Clausen's theorem and Sun's congruence $(24)$, which is itself a consequence of Bernoulli's formula and von Staudt-Clausen's theorem. The a priori non trivial equivalence between our formula and Sun's formula dating from $2000$ got established in the introduction of the current paper.

\section{Concluding words}

Because the introduction was quite comparative and interactive, by times intertwining, between Glaisher's pioneering and subtle work, Sun's technical and highly beautiful work a hundred years later and our current work another twenty years later, we felt the need to gather a final conclusion regarding the different approaches. Thus, the forthcoming lines echo very much what was said in the lengthy introduction and draw a short synthesis.

First we discuss about Wilson's theorem. It appears to us that our approach is the most natural one and the most straightforward one for proving Glaisher's original result from $1900$ on Wilson's theorem modulo $p^2$. Indeed, Glaisher's own proof for the factorial uses the knowledge of the Stirling numbers modulo $p^2$ which he computes by quite original means. Moreoever, Glaisher's method does not generalize to the modulus $p^3$. Sun's method leads to the next modulus $p^3$, but it involves a lot of different works like finding the generalized harmonic numbers modulo $p^3$. Also, his method does not generalize further, to the modulus $p^4$ for instance, while ours does.

Now we discuss about the three works involving the Stirling numbers and the generalized harmonic numbers, and their interconnection. \\
Sun's way based on Newton's formulas for computing the Stirling numbers modulo $p^2$ is by far the most elegant and efficient. Once this is achieved, our method allows to relate modulo $p^2$ the Stirling numbers as a difference of a sum of powers and of a generalized harmonic number. We use Bayat's result (only modulo $p$) which is a generalization of Wolstenholme's theorem to generalized harmonic numbers. Thus, we deduce the generalized harmonic numbers modulo $p^2$ by a simpler method than Sun's. However, our method this time does not generalize to the modulus $p^3$ whereas Sun's sophisticated method does. In fact his method exposed in \cite{SU2} even generalizes to the modulus $p^4$ (the result gets provided in his Remark $5.1$ without details of his calculation). Also, our method does not deal with $H_{p-1,p-1}$ whereas Sun treats that case as well.
Our method allows to improve Bayat's result (or Hardy and Wright's) modulo $p^2$ since we show that $H_{p-1,p-2}=0\;\text{mod}\,p^2$ and that case is not contained in the works of the latter authors. Recall that in order to obtain our result we only use Bayat's result modulo $p$. \\
Going back to the Stirling numbers, Sun's method does not generalize to the modulus $p^3$, unfortunately. In fact, it does, but only for the even indices (with respect to our notations). This was done by us in $\S\,1$, see Eq. $(27)$ using his strategy modulo $p^2$. Glaisher also obtains the Stirling numbers modulo $p^3$, again limited to the even indices. In a forthcoming paper, we will deal with all the indices, even and odd, following the method presented in this paper. As part of our work, we will resolve some congruences concerning some convolutions involving Bernoulli numbers. Glaisher's work for finding the Stirling numbers modulo $p^3$ for even indices leads him to knowing the generalized harmonic numbers modulo $p^3$ for odd indices. Namely, Glaisher succeeds to relate the Stirling numbers and the harmonic numbers modulo $p^2$ (resp $p^3$) for even (resp odd) indices (with respect this time to his own notations). His proof involves a smart use of Newton's formulas and knowing the $p^2$ divisibility property of the Stirling numbers with even indices. He knew these divisibility properties from his original work on the Stirling numbers modulo $p^2$ which had led to his $1900$ formula for the factorial modulo $p^2$. We note once again like in $\S\,2$ that an alternative to Sun's way for finding the Stirling numbers modulo $p^2$ is our method joint with Bayat's result modulo $p^2$ for determining the Stirling numbers with even indices modulo $p^2$. Then the odd indices can be resolved by using Glaisher's relationship and our relationship from Theorem $3$. The least efficient way for determining the Stirling numbers modulo $p^2$ consists of using our Theorem $3$ and Sun's congruences for generalized harmonic numbers modulo $p^2$.

In brief, Wilson's theorem should be worked out using our $p$-adic method. The Stirling numbers modulo $p^2$ should be determined following Sun.
Then the generalized harmonic numbers modulo $p^2$ can be derived by our Theorem $3$. For odd indices, Glaisher's way is the best to determine the generalized harmonic numbers modulo $p^3$. For even indices, Sun's heavy machinery can't be bypassed. The Stirling numbers modulo $p^3$ with even indices can be obtained by our method built on Sun's method, or can be determined from the Stirling numbers modulo $p^2$ with odd indices by using Glaisher's method. The Stirling numbers with odd indices modulo $p^3$ would benefit from further investigations. \\

The research related to the current paper has drawn the attention of mathematicians over several centuries and is still ongoing research.
A sign also of the vibrant interest in the topics presented here is the fact that many authors have been working on quantum analogs.
In the literature, people have for instance studied quantum harmonic numbers
$$H_{p-1}(q)=\sum_{k=1}^{p-1}\frac{1}{[k]_q}\;\;\text{and}\;|q|<1,$$
with $[k]_q$ the so-called quantum integer defined by $$[k]_q=\frac{1-q^k}{1-q}=1+q+\dots+q^{k-1}$$
In $1999$, G.E. Andrews shows in \cite{AN} the interesting congruence
\begin{Result} (Andrews 1999 \cite{AN})
\begin{equation}
H_{p-1}(q)\equiv\,\frac{p-1}{2}\,(1-q)\;\;\text{mod}\,[p]_q
\end{equation}
\end{Result}
\noindent His work gets enriched modulo the squared quantum integer by L.L. Shi and H. Pan in $2007$ \cite{SP}. Shi and Pan's result is even a generalization of Wolstenholme's theorem at the next higher power of quantum integers, namely
\begin{Result} (Shi and Pan 2007 \cite{SP})
$$\left|\begin{array}{ccccc}
H_{p-1}(q)&\equiv&\frac{p-1}{2}(1-q)+\frac{p^2-1}{24}(1-q)^2\,[p]_q&\text{mod}&[p]_q^2\\
&&&&\\
H_{p-1,2}(q)&\equiv&-\frac{(p-1)(p-5)}{12}\,(1-q)^2&\text{mod}&[p]_q
\end{array}\right.$$
\end{Result}
\noindent And in fact, in \cite{DI}, Dilcher even tackles congruences with generalized quantum harmonic numbers. He defines two key determinants both involving binomial coefficients.
\begin{Definition} (Dilcher's determinants) Let
$$
D_k(p)=\left|\begin{array}{ccccc}
\bn{p+1}{2}&p&0&\dots&0\\
\bn{p+1}{3}&\bn{p+1}{2}&p&\dots&0\\
\vdots&\vdots&\vdots&\ddots&\vdots\\
\bn{p+1}{k}&\bn{p+1}{k-1}&\bn{p+1}{k-2}&\dots&p\\
\bn{p+1}{k+1}&\bn{p+1}{k}&\bn{p+1}{k-1}&\dots&\bn{p+1}{2}
\end{array}\right|$$
and $$\tilde{D}_k(p)=\left|\begin{array}{ccccc}
\bn{p}{2}&\bn{p}{1}&0&\dots&0\\
\bn{p}{3}&\bn{p}{2}&\bn{p}{1}&\dots&0\\
\vdots&\vdots&\vdots&\ddots&\vdots\\
\bn{p}{k}&\bn{p}{k-1}&\bn{p}{k-2}&\dots&\bn{p}{1}\\
\bn{p}{k+1}&\bn{p}{k}&\bn{p}{k-1}&\dots&\bn{p}{2}
\end{array}\right|$$
\end{Definition}
\begin{Notation}
$$H_{p-1,k}(q):=\sum_{j=1}^{p-1}\frac{1}{[j]_q^k}\;\;\text{and}\;\;\tilde{H}_{p-1,k}(q):=\sum_{j=1}^{p-1}\frac{q^j}{[j]_q^k},$$
\end{Notation}
\noindent From Dilcher, the second sum is more natural to deal with. It is the reason why Dilcher introduces it. Dilcher shows that
\begin{Result} (Dilcher 2008 \cite{DI})
$$\begin{array}{l}
H_{p-1,k}(q)\equiv \frac{(-1)^{k-1}}{p^k}\,D_k(-p)(1-q)^k\;(mod\,[p]_q)\\\\
\tilde{H}_{p-1,k}(q)\equiv-\frac{1}{p^k}\,\tilde{D}_k(p)(1-q)^k\;(mod\,[p]_q)
\end{array}$$
\end{Result}
\noindent Even more recently in $2015$, researchers at the Fields Institute deal with q-analogues for congruences involving multiple harmonic sums,
after dealing themselves with multiple harmonic sums in a non quantum setting \cite{PI1}. In \cite{PI2}, the authors introduce multiple q-harmonic sums which they define as $$H_{p-1}(\mathbf{s},\mathbf{t})=\sum_{1\leq k_1<\dots<k_l\leq p-1}\frac{q^{k_1t_1+\dots+k_kt_l}}{[k_1]_q^{s_1}\dots [k_l]_q^{s_l}}$$
for two $l$-tuples of non-negative integers $\mathbf{s}:=(s_1,\dots,s_l)$ and $\mathbf{t}:=(t_1,\dots,t_l)$.

Binomial congruences also have their quantum analogs, including Wolstenholme's congruence. Interest in binomial congruences goes back to K\"ummer and Lucas.\\
K\"ummer shows in $1852$ that if $p^r$ is the highest power of $p$ dividing $\binom{n}{m}$, then $r$ equals the number of carries when adding $m$ and $n-m$ in base $p$ arithmetic. \\
Lucas in $1878$ shows that the binomial coefficients modulo $p$ are linked to the respective expansions of the two integers in base $p$.
\begin{Thm}(Lucas $1878$)\\
If $$\begin{array}{l} n=n_0+n_1p+\dots +n_sp^s\\m=m_0+m_1p+\dots +m_sp^s\end{array}$$
then, $$\binom{n}{m}\equiv\,\prod_{i=0}^s\binom{n_i}{m_i}\;\text{mod}\,p$$
\end{Thm}
\noindent A quantum analogue of the binomial Wolstenholme theorem does exist since $2011$, provided by A. Straub \cite{ST}. \\
An equivalent form for Wolstenholme's theorem is the following.
$$\binom{2p}{p}\equiv 2\;\text{mod}\,p^3$$
Congruences involving binomial coefficients $\binom{np}{mp}$ for $p$ prime have been studied by mathematicians for decades. The most recent result that is now known is due to Helou and Terjanian \cite{HT}. They showed in $2008$ that for $p$ prime with $p\geq 5$ and integers $n$ and $m$ such that $0<m\leq n$, we have
$$\binom{np}{mp}\equiv\,\binom{n}{m}\;\text{mod}\,p^s$$
with $s$ the highest power of $p$ dividing $p^3\,m(n-m)\binom{n}{m}$.\\
Quantum wise, Clark was the first to tackle quantum binomial congruences. \\He showed in $1995$ \cite{CA} that
$$\binom{np}{mp}_q\equiv\,\binom{n}{m}_{q^{p^2}}\;\text{mod}\,[p]_q^2,$$
where the binomial coefficients are defined with the quantum factorials. \\
Fifteen years later, Straub improves Clark's congruence by showing under the same conditions as before that
$$\binom{np}{mp}\equiv\binom{n}{m}_{q^{p^2}}-\binom{n}{m+1}\binom{m+1}{2}\frac{p^2-1}{12}(q^p-1)^2\;\text{mod}[p]_q^3$$
The case $n=2$ and $m=1$ provides a q-analogue for the binomial Wolstenholme theorem, namely
\begin{Thm} (Straub $2011$, q-analogue for Wolstenholme's theorem)
\begin{equation}\binom{2p}{p}_q\equiv [2]_{q^{p^2}}-\frac{p^2-1}{12}\,(q^p-1)^2\;\text{mod}\,[p]_q^3\end{equation}
\end{Thm}
\noindent Note that Andrews \cite{AN} had already contributed modulo the squared power by showing that
$$\binom{2p-1}{p-1}_q\equiv\,q^{\frac{p(p-1)}{2}}\;\text{mod}\,[p]^2$$
It is to expect that there will be many more generalizations of these results. Most importantly, in what setting could we use these generalizations ? Time shall also provide an answer to that question.

On a final note, going back to the recently honored mathematician from the title, it is quite pleasant to mention that Faulhaber has a quantum analogue for his formula that is due to V.J.W. Guo and J. Zeng, see \cite{GZ}.\\
Also, classically, it would be nice to understand his coefficients better $p$-adically. We think that the methods of this paper can contribute to it.
\\

\textit{Email address:} \textit{clairelevaillant@yahoo.fr}

\end{document}